\documentclass[journal]{IEEEtran}
\ifCLASSINFOpdf
   \usepackage[pdftex]{graphicx}
\else
\fi
\usepackage{bm,algorithmicx,algpseudocode,threeparttable,booktabs}
\usepackage{amsmath,color}
\usepackage{cases}
\usepackage{stmaryrd}
\usepackage{amsfonts}
\usepackage{graphicx,times,psfig,amsmath}
\usepackage{cite}
\usepackage{hyperref}
\newtheorem{remark}{Remark}
\DeclareMathOperator*{\argmin}{arg\,min}

\begin{document}
%
\title{A Statistical Study on Parameter Selection of Operators in Continuous State Transition Algorithm}
%
%
%

\author{Xiaojun~Zhou, Chuhua Yang ~\IEEEmembership{Member,~IEEE}, Weihua Gui ~\IEEEmembership{Member,~IEEE,}
\thanks{X. Zhou, C. Yang and W. Gui are with the School
of Information Science and Engineering, Central South University, Changsha, Hunan, China e-mail: (michael.x.zhou@csu.edu.cn, ychh@csu.edu.cn, gwh@csu.edu.cn).}
\thanks{Manuscript received XX, 2017; revised XX, 2018.}}

\maketitle

\begin{abstract}
State transition algorithm (STA) has been emerging as a novel metaheuristic method for global optimization in recent few years.
In our previous study, the parameter of transformation operator in continuous STA is kept constant or decreasing itself in a periodical way.
In this paper, the optimal parameter selection of the STA is taken in consideration. Firstly, a statistical study with four benchmark two-dimensional  functions is conducted to show how these parameters affect the search ability of the STA. Based on the experience gained from the statistical study, then, a new continuous STA with optimal parameters strategy is proposed to accelerate its search process. The proposed STA is successfully applied to twelve benchmarks with 20, 30 and 50 dimensional space. Comparison with other metaheuristics has also demonstrated the effectiveness of the proposed method.
\end{abstract}


%
\IEEEpeerreviewmaketitle

\section{Introduction}
\IEEEPARstart{S}{tate-transition-algorithm} (STA) \cite{zhou2012state,zhou2016discrete} is a recently emerging
metaheuristic method for global optimization and has found applications in nonlinear system identification and control \cite{zhou2014nonlinear},
water distribution networks configuration \cite{zhou2016optimal}, sensor network localization \cite{zhou2018dynamic}, PID controller design \cite{zhang2016fractional,saravanakumar2015tuning}, overlapping peaks resolution \cite{wang2016state}, image segmentation \cite{han2017new}, wind power prediction \cite{wang2016new},
dynamic optimization \cite{han2017dynamic, huang2017dynamic}, bi-level optimization \cite{huang2018novel}, modeling and control of complex industrial processes
\cite{wang2016optimization,xie2016new,xie2018data,yang2017optimal,zhang2018fractional,zhou2018setpoint}, etc.
In STA, a solution to
an optimization problem is considered as a state, and an update of a solution can be
regarded as a state transition. Unlike the population-based evolutionary algorithms \cite{wang2018cooperative,wang2018global, zhang2018dual}, the standard STA is an individual-based optimization method.
Based on an incumbent best solution, a neighborhood with special characteristics will be formed automatically
when using certain state transformation operator.
A variety of state transformation operators, for example,
rotation, translation, expansion, and axesion in continuous STA, or swap, shift, symmetry and
substitute in discrete STA,
are designed purposely for both global and local search.
On the basis of the neighborhood,
then, a sampling technique is used to generate a candidate set, and
the next best solution is updated by using a selection technique based on previous best solution and
the candidate set. This process is repeated using state transformation operators alternatively until some terminal conditions are satisfied.

In this paper, the continuous state transition algorithm is studied.
As aforementioned, in continuous STA, there are four state transformation operators, and each transformation operator has certain geometric significance, i.e., the neighborhood formed by each
transformation operator has certain geometric characteristic. To be more specific,
the rotation transformation has the functionality to search in a hypersphere with the maximal radius $\alpha$, called rotation factor;  the translation transformation has the functionality to search along
a line with the maximal length $\beta$, called translation factor;
the expansion transformation has the functionality to search in a broader space controlled by the expansion factor $\gamma$; and the axesion
transformation is designed to strengthen single-dimensional search regulated by the axesion factor $\delta$.
In our previous studies, the rotation factor is exponentially decreasing from a maximum value to a minimum value in a
periodic way, and other transformation factors are kept constant at one \cite{zhou2012state}.
To gain a better exploitation ability, all state transformation factors are exponentially decreasing from a maximum value to a minimum value in a
periodic way in \cite{zhou2018dynamic}.

As is known to us, there exist several parameters in metaheuristic methods and parameter selection plays a significant role in their performance.
For instance, crossover and mutation probability in genetic algorithms (GAs) \cite{eiben2011parameter}, inertia weight and acceleration
factors in particle swarm optimization (PSO) \cite{sun2012quantum,zhang2014parameter},
amplification factor and crossover rate in differential evolution (DE) \cite{mallipeddi2011differential,sarker2014differential,wang2011differential}, and
neighborhood radius in artificial bee colony (ABC) \cite{karaboga2014quick}. In general, the parameter setting can be summarized to two types: parameter tuning and parameter control. The former is to find good parameter values before running these algorithms, and they remain fixed during the run. On the contrary, the later is to update parameter values in the process, and the types of
update mechanisms can be deterministic, adaptive, or self-adaptive (for details, please refer to \cite{eiben1999parameter, de2007parameter, karafotias2015parameter}).

To gain a better understanding of the parameters of transformation operators in continuous STA affecting its performance, the parameter selection in continuous STA is focused in this study. With four commonly used benchmark functions as cases, several properties of the operator parameters are observed from a statistical study.
With the gained experience from the statistical results, a new continuous STA with optimal operator parameter selection strategy is proposed, and the proposed STA is successfully applied to other benchmarks with higher dimensions.

The remainder of this paper is organized as follows.
In Section II, the standard continuous STA are described.
Section III gives a statistical study to show how the operator parameters in continuous STA affecting its performance .
The proposed STA with optimal operator parameter selection strategy is given in Section IV.
In Section V, experimental results are given to testify the effectiveness of the proposed STA.
Finally, conclusion is drawn in Section VI.

\section{Standard continuous state transition algorithm}
Consider the following continuous optimization problem with simple constraints:
\begin{eqnarray}
\min_{\bm x \in \Omega} f(\bm x)
\end{eqnarray}
where $\Omega \subseteq \mathbb{R}^n$ is a closed and compact set, which is usually composed of lower and upper bounds of $\bm x$, \textit{i.e.}, $\Omega = \{\bm x \in \mathbb{R}^n| \underline{x_i} \leq x_i \leq \overline{x_i}, i = 1, \cdots, n\}$.

In classical iterative methods for numerical optimization, a new candidate is generated based on a previous solution by using different optimization operators. In a state transition way, a solution can be regarded as a state, and an update of a solution can be considered as a state transition. On the basis of state space representation, the unified form of generation of solution in state transition algorithm can be described as follows:
\begin{equation}
\left \{ \begin{array}{ll}
\bm s_{k+1}= A_{k} \bm s_{k} + B_{k} \bm u_{k}\\
y_{k+1}= f(\bm s_{k+1})
\end{array} \right.,
\end{equation}
where $\bm s_{k}$ and $\bm s_{k+1}$  stand for a current state and the next state respectively, corresponding to solutions of the optimization problem; $\bm u_{k}$ is a function of $\bm s_{k}$ and historical states; $y_{k}$ is the fitness value at $\bm s_k$;
$A_{k}$ and
$B_{k}$ are state transition matrices, which can be considered as transformation operators;
 $f$ is the objective function or fitness function.

\subsection{State transition operators}
Using state space representation and state transformation for reference, four special
state transformation operators are designed to generate candidate solutions for an optimization problem \cite{zhou2011initial, zhou2012state}.\\
(1) Rotation transformation
\begin{equation}
\bm s_{k+1}= \bm s_{k}+\alpha \frac{1}{n \|\bm s_{k}\|_{2}} R_{r} \bm s_{k},
\end{equation}
where $\alpha$ is a positive constant, called the rotation factor;
$R_{r}$ $\in$ $\mathbb{R}^{n\times n}$, is a random matrix with its entries being uniformly distributed random variables defined on the interval [-1, 1],
and $\|\cdot\|_{2}$ is the L2-norm (or Euclidean norm) of a vector. This rotation transformation
has the functionality to search in a hypersphere with the maximal radius $\alpha$, which has been testified. The rotation transformation is designed for local search and can be used to guarantee local optimality and manipulate solution accuracy. \\
(2) Translation transformation\\
\begin{equation}
\bm s_{k+1} = \bm s_{k}+  \beta  R_{t}  \frac{\bm s_{k}- \bm s_{k-1}}{\|\bm s_{k}- \bm s_{k-1}\|_{2}},
\end{equation}
where $\beta$ is a positive constant, called the translation factor; $R_{t}$ $\in \mathbb{R}$ is a uniformly distributed random variable defined on the interval [0,1].
It is not difficult to understand that the translation transformation has the functionality to search along a line from $s_{k-1}$ to $s_{k}$ at the starting point $s_{k}$ with maximum length $\beta$. The translation operator is actually a line search, and it can be considered as a heuristic operator since there exists a possible better solution along the line if $s_k$ is better than $s_{k-1}$.
\\
(3) Expansion transformation\\
\begin{equation}
\bm s_{k+1} = \bm s_{k}+  \gamma  R_{e} \bm s_{k},
\end{equation}
where $\gamma$ is a positive constant, called the expansion factor; $R_{e} \in \mathbb{R}^{n \times n}$ is a random diagonal
matrix with its entries obeying the Gaussian distribution (or  normal distribution). In the standard STA, the mean equals zero and standard deviation equals one, i.e., the standard normal distribution is used.
The expansion transformation has the functionality to search in the whole space in probability, and it is designed for global search.
\\
(4) Axesion transformation\\
\begin{equation}
\bm s_{k+1} = \bm s_{k}+  \delta  R_{a}  \bm s_{k}\\
\end{equation}
where $\delta$ is a positive constant, called the axesion factor; $R_{a}$ $\in \mathbb{R}^{n \times n}$ is a random diagonal matrix with its entries obeying the Gaussian distribution and only one random position having nonzero value.
The axesion
transformation is designed to search along the axes, aiming to strengthen single-dimensional search \cite{zhou2011new}.

\subsection{A sampling technique}
The idea of sampling incorporated in continuous STA was firstly illustrated in \cite{zhou2016comparative}.
It is found that for a given solution, a neighborhood will be automatically formed.
To avoid enumerating all possible candidate solutions, representative samples can be used
to reflect the characteristics of the neighborhood. Taking the rotation transformation for example, when independently executing the rotation operator for \textit{SE} times,  a total number of
\textit{SE} samples are generated in pseudocode as follows
\begin{algorithmic}[1]
\For{$i\gets 1$, \textit{SE}}
\State $\mathrm{State}(:,i) \gets {\mathrm{Best}} +\alpha \frac{1}{n \| {\mathrm{Best}} \|_{2}} R_{r}  {\mathrm{Best}} $
\EndFor
\end{algorithmic}
where $Best$ is the incumbent best solution, and \textit{SE} samples are stored in
the matrix $\mathrm{State}$.

\subsection{An update strategy}
As mentioned above, based on the incumbent best solution, a total number of
\textit{SE} candidate solutions are generated, but it should be noted that these candidate solutions do not always belong to the domain $\Omega$.
To address this issue, these samples are projected into
$\Omega$ through
\begin{equation}
x_i =
\left\{ \begin{aligned}
&\overline{x_i},  \;\;\; \mathrm{if}\; x_i > \overline{x_i} \\
&\underline{x_i},  \;\;\; \mathrm{if}\; x_i < \underline{x_i}\\
&x_i,  \;\;\; \mathrm{otherwise}
\end{aligned} \right.
\end{equation}

As a result, the candidate solutions can be guaranteed to be always feasible.
Next, a new best solution is selected from the candidate set by virtue of the fitness function, denoted as $newBest$.
Finally, an update strategy based on greedy criterion is used to update the incumbent best as shown below
\begin{equation}
Best =
\left\{ \begin{aligned}
&\mathrm{newBest},  \;\;\;\mathrm{if}\; f(\mathrm{newBest}) < f(\mathrm{Best}) \\
&\mathrm{Best},\;\;\;\;\;\;    \;\;\;\mathrm{otherwise}
\end{aligned} \right.
\end{equation}

\subsection{Algorithm procedure of the standard continuous STA}
With the state transformation operators for both local and global search, sampling technique for time-saving and update strategy for convergence, the standard continuous STA can be described by the following pseudocodes
\begin{algorithmic}[1]
\State State $\gets$ initialization(SE,$\Omega$)
\State Best $\gets$ fitness(funfcn,State)
\Repeat
    \If{$\alpha < \alpha_{\min}$}
    \State {$\alpha \gets \alpha_{\max}$}
    \EndIf
    \State {Best $\gets$ expansion(funfcn,Best,SE,$\beta$,$\gamma$)}
    \State {Best $\gets$ rotation(funfcn,Best,SE,$\beta$,$\alpha$)}
    \State {Best $\gets$ axesion(funfcn,Best,SE,$\beta$,$\delta$)}
    \State {$\alpha \gets \frac{\alpha}{\textit{fc}}$}
\Until{the specified termination criterion is met}
\end{algorithmic}

\indent As for detailed explanations, rotation$(\cdot)$ in above pseudocode is given for illustration purposes as follows
\begin{algorithmic}[1]
\State{oldBest $\gets$ Best}
\State{fBest $\gets$ feval(funfcn,oldBest)}
\State{State $\gets$ op\_rotate(Best,SE,$\alpha$)}
\State{[newBest,fnewBest] $\gets$ fitness(funfcn,State)}
\If{fnewBest $<$ fBest}
    \State{fBest $\gets$ fnewBest}
    \State{Best $\gets$ newBest}
    \State{State $\gets$ op\_translate(oldBest,newBest,SE,$\beta$)}
    \State{[newBest,fnewBest] $\gets$ fitness(funfcn,State)}
    \If{fnewBest $<$ fBest}
        \State{fBest $\gets$ fnewBest}
        \State{Best $\gets$ newBest}
    \EndIf
\EndIf
\end{algorithmic}

As shown in the above pseudocodes, initialization $(\cdot)$ is used to make sure the
initial solution is in the range $\Omega$. The rotation factor $\alpha$ is decreasing periodically from a maximum
value $\alpha_{\max}$ to a minimum value $\alpha_{\min}$ in an
exponential way with base \textit{fc}, which is called lessening coefficient.
op\_rotate$(\cdot)$ and op\_translate$(\cdot)$ represent the implementations of proposed sampling technique for rotation and
translation operators, respectively, and fitness$(\cdot)$ represents the implementation of selecting the new best solution
from \textit{SE} samples. It should be emphasized that the translation operator is only executed
when a solution better than the incumbent best solution can be found in the \textit{SE} samples
from rotation, expansion or axesion transformation.
In the standard continuous STA, the parameter settings are given as follows:
$\alpha_{\max} = 1, \alpha_{\min} = 1e$-4, $\beta = 1, \gamma = 1, \delta = 1$, $\textit{SE} = 30, \textit{fc} = 2$.

\section{Statistical study of the state transformation factors}
As described in Section II, in the standard continuous STA,  the state transformation factors like
expansion factor $\gamma$, axesion factor $\delta$ are kept constant, and
rotation factor $\alpha$ is decreasing periodically from a maximum
value $\alpha_{\max}$ to a minimum value $\alpha_{\min}$ in an
exponential way.
In order to select the values of these parameters in a more effective manner,
a statistical study of the state transformation factors is carried out to investigate the effect of
parameter selection on the performance of state transition operators.

Four well-known benchmark functions are listed below:

(1) Spherical function
\[
f_1(\bm x)= \sum_{i=1}^n x_{i}^{2},
\]
where the global optimum $\bm x^{*} = (0, \cdots, 0)$ and $f(\bm x^{*}) = 0$,
$-100 \leq x_i \leq 100, i= 1, \cdots, n$.

(2) Rosenbrock function
\[
f_2(\bm x)= \sum_{i=1}^n (100(x_{i+1}-x_{i}^2)^2 + (x_{i}-1)^2),
\]
where the global optimum $\bm x^{*} = (1, \cdots, 1)$ and $f(\bm x^{*}) = 0$,
$-30 \leq x_i \leq 30, i= 1, \cdots, n$.

(3) Rastrigin function
\[
f_3(\bm x) =  \sum_{i=1}^n(x_{i}^{2}-10\cos(2 \pi x_{i})+10),
\]
where the global optimum $\bm x^{*} = (0, \cdots, 0)$ and $f(\bm x^{*}) = 0$,
$-5.12 \leq x_i \leq 5.12, i= 1, \cdots, n$.

(4) Griewank function
\[
f_4(\bm x) = \frac{1}{4000} \sum_{i=1}^n x_{i}^{2}-  \prod_{i}^{n} \cos|\frac{x_{i}}{\sqrt{i}}| + 1,
\]
where the global optimum $\bm x^{*} = (0, \cdots, 0)$ and $f(\bm x^{*}) = 0$,
$-600 \leq x_i \leq 600, i= 1, \cdots, n$.

For a given solution $Best_0$, three state transition operators (rotation, expansion and axesion) are performed respectively
for \textit{SE} times (yielding \textit{SE} samples) independently
on each benchmark function
using different values of state transformation factors.
To be more specific, there are five groups of given solutions, \textit{i.e.},
$Best_0 = (0.01,0.01), (0.1,0.1), (0.5,0.5), (0.9, 0.9)$,  $(0.99, 0.99)$;
the total number of samples is set at \textit{SE} = 1e6; and the value of state transformation operators
is chosen from the set $\Omega = \{$1, 1e-1, 1e-2, 1e-3, 1e-4, 1e-5, 1e-6, 1e-7, 1e-8$\}$.

To evaluate the influence of the parameter selection on the performance of state transition operators,
the following two indexes are introduced:
\begin{eqnarray}
\rho_s &=& \frac{N_s}{\mathrm{\textit{SE}}} \\
\rho_d &=& \frac{|ave - fBest_0|}{|fBest_0|}
\end{eqnarray}
where $\rho_s$ and $\rho_d$ are called success rate and descent rate, respectively. $N_s$ is the number
of samples whose objective function values are smaller than that of the $Best_0$.
$ave$ is the average function value of the $N_s$ samples, and
$fBest_0$ represents the
function value for $Best_0$,

The statistical results for different values of state transformation factors can be found from Table I to Table IV.

\begingroup
\renewcommand{\baselinestretch}{1}%
\begin{table*}[htbp!]
\centering
\begin{threeparttable}[b]
\renewcommand{\arraystretch}{1.3}
\caption{Statistical results of success rate and descent rate for the rotation transformation (Spherical problem)}
\scriptsize
\label{DSTA_benchmarks}
\begin{tabular}{{p{1.2cm}|p{0.6cm}|p{1.25cm}|p{1.25cm}|p{1.25cm}|p{1.25cm}|p{1.25cm}|p{1.25cm}|p{1.25cm}|p{1.25cm}|p{1.25cm}}}
\hline
$\mathrm{Best}_0$ & Index & $\alpha$ = 1 & $\alpha$ = 0.1 & $\alpha$ = 0.01 & $\alpha$ = 1e-3 & $\alpha$ = 1e-4 & $\alpha$ = 1e-5 & $\alpha$ = 1e-6 & $\alpha$ = 1e-7 & $\alpha$ = 1e-8\\
\hline
$(0.01,0.01)$ & $\rho_s$ &  0.0012 &   0.0919 &   0.4550  &  0.4953  &  0.5002  &  0.5000 &   0.4991 &   0.5001 &   0.4999  \\
              & $\rho_d$ & 5.0448e-1 & 5.0347e-1 & 2.7785e-1 & 3.2465e-2  &3.2970e-3 & 3.2960e-4 & 3.2990e-5  & 3.2923e-6 & 3.3056e-7\\
\hline
$(0.1,0.1)$  & $\rho_s$ &  0.0921  &  0.4555 &   0.4956 &   0.4990  &  0.4996  &  0.4993  &  0.4996 &   0.5001   & 0.5012\\
             & $\rho_d$ & 5.0313e-1&  2.7801e-1  &3.2470e-2 & 3.2949e-3 & 3.3048e-4 & 3.3026e-5&  3.2996e-6 & 3.2944e-7 & 3.2962e-8\\
\hline
$(0.5,0.5)$  & $\rho_s$ & 0.4144  &  0.4899  &  0.4992  &  0.5000  &  0.4999 &   0.5006 &   0.4989 &   0.5008 &   0.5003
\\
             & $\rho_d$ & 4.4908e-1 & 6.3938e-2 & 6.5821e-3 & 6.5914e-4 & 6.5918e-5 & 6.6006e-6&  6.5986e-7 & 6.5950e-8 & 6.6128e-9
             \\
\hline
$(0.9,0.9)$  & $\rho_s$ & 0.4498  &  0.4947  &  0.4994 &   0.5002  &  0.5000  &  0.4994 &   0.4995  &  0.4996 &   0.4989\\
             & $\rho_d$ & 3.0238e-1 & 3.6061e-2 & 3.6623e-3 & 3.6555e-4 & 3.6682e-5 & 3.6639e-6 & 3.6656e-7  & 3.6734e-8 & 3.6637e-9\\
\hline
$(0.99,0.99)$  & $\rho_s$ & 0.4543  &  0.4955 &   0.4986 &   0.5000 &   0.4995 &   0.4999 &   0.5001 &   0.5008    & 0.4999\\
            & $\rho_d$  & 2.7968e-1 & 3.2731e-2 & 3.3249e-3 & 3.3378e-4 & 3.3307e-5&  3.3289e-6&  3.3337e-7  & 3.3349e-8 & 3.3323e-9\\
\hline
\end{tabular}
\end{threeparttable}
\end{table*}
\endgroup

\begingroup
\renewcommand{\baselinestretch}{1}%
\begin{table*}[htbp!]
\centering
\begin{threeparttable}[b]
\renewcommand{\arraystretch}{1.3}
\caption{Statistical results of success rate and descent rate for the rotation transformation (Rosenbrock problem)}
\scriptsize
\label{DSTA_benchmarks}
\begin{tabular}{{p{1.2cm}|p{0.6cm}|p{1.25cm}|p{1.25cm}|p{1.25cm}|p{1.25cm}|p{1.25cm}|p{1.25cm}|p{1.25cm}|p{1.25cm}|p{1.25cm}}}
\hline
$\mathrm{Best}_0$ & Index & $\alpha$ = 1 & $\alpha$ = 0.1 & $\alpha$ = 0.01 & $\alpha$ = 1e-3 & $\alpha$ = 1e-4 & $\alpha$ = 1e-5 & $\alpha$ = 1e-6 & $\alpha$ = 1e-7 & $\alpha$ = 1e-8\\
\hline
$(0.01,0.01)$ & $\rho_s$ &  0.0700  &  0.2665 &   0.4766 &   0.4985  &  0.5001 &   0.4997 &   0.5005  &  0.4991 &   0.4999\\
              & $\rho_d$ & 2.9603e-1 & 4.3729e-2 & 6.1371e-3 & 6.6152e-4 & 6.6718e-5 & 6.6619e-6 & 6.6708e-7 & 6.6658e-8 & 6.6683e-9\\
\hline
$(0.1,0.1)$  & $\rho_s$ & 0.1891  &  0.5082 &   0.5019 &   0.4996  &  0.5003  &  0.5000 &   0.5005  &  0.4993  &  0.4999\\
             & $\rho_d$ & 3.9041e-1 & 2.2187e-1 & 2.6801e-2 & 2.7281e-3 & 2.7272e-4 & 2.7264e-5 & 2.7332e-6  & 2.7328e-7 & 2.7319e-8\\
\hline
$(0.5,0.5)$  & $\rho_s$ & 0.3884  &  0.5049   & 0.5005  &  0.5003  &  0.5003 &   0.5004 &   0.4996  &  0.5010  &  0.4997
\\
             & $\rho_d$ & 6.4979e-1 & 2.3369e-1 & 2.5464e-2 & 2.5633e-3 & 2.5620e-4 & 2.5615e-5 & 2.5619e-6 & 2.5595e-7 & 2.5623e-8\\
\hline
$(0.9,0.9)$  & $\rho_s$ & 0.1121  &  0.5006  &  0.4999  &  0.5004  &  0.5001  &  0.4999  &  0.5005 &   0.4994  &  0.4996\\
             & $\rho_d$ & 6.4261e-1 & 6.3423e-1 & 1.0225e-1 & 1.0622e-2 & 1.0642e-3 & 1.0650e-4 & 1.0645e-5  & 1.0649e-6 & 1.0633e-7\\
\hline
$(0.99,0.99)$  & $\rho_s$ & 0.0052 &   0.1146 &  0.5004 &   0.5003  &  0.4996  &  0.4992 &   0.5003 &   0.4994   & 0.5003\\
            & $\rho_d$  & 5.1273e-1 & 6.5350e-1 & 6.2945e-1 & 1.0034e-1 & 1.0419e-2 & 1.0446e-3 & 1.0448e-4 &  1.0457e-5 & 1.0441e-6\\
\hline
\end{tabular}
\end{threeparttable}
\end{table*}
\endgroup

\begingroup
\renewcommand{\baselinestretch}{1}%
\begin{table*}[htbp!]
\centering
\begin{threeparttable}[b]
\renewcommand{\arraystretch}{1.3}
\caption{Statistical results of success rate and descent rate for the rotation transformation (Rastrigin problem)}
\scriptsize
\label{DSTA_benchmarks}
\begin{tabular}{{p{1.2cm}|p{0.6cm}|p{1.25cm}|p{1.25cm}|p{1.25cm}|p{1.25cm}|p{1.25cm}|p{1.25cm}|p{1.25cm}|p{1.25cm}|p{1.25cm}}}
\hline
$\mathrm{Best}_0$ & Index & $\alpha$ = 1 & $\alpha$ = 0.1 & $\alpha$ = 0.01 & $\alpha$ = 1e-3 & $\alpha$ = 1e-4 & $\alpha$ = 1e-5 & $\alpha$ = 1e-6 & $\alpha$ = 1e-7 & $\alpha$ = 1e-8\\
\hline
$(0.01,0.01)$ & $\rho_s$ &  0.0012 &   0.0918 &   0.4552  &  0.4948 &   0.4993 &   0.5000 &   0.5003 &   0.5000 &   0.4995\\
              & $\rho_d$ & 5.0017e-1 &  5.0278e-1 &  2.7796e-1 &  3.2470e-2 &  3.2959e-3 &  3.3012e-4  & 3.2908e-5 &  3.3004e-6 &  3.3061e-7\\
\hline
$(0.1,0.1)$  & $\rho_s$ & 0.0932  &  0.4606  &  0.4960  &  0.4992  &  0.5000  &  0.4995  &  0.4997 &   0.4996  &  0.4998\\
             & $\rho_d$ & 4.9612e-1 &  2.7255e-1 &  3.1446e-2 &  3.1888e-3 &  3.1848e-4 &  3.1876e-5  & 3.1924e-6 &  3.1915e-7 &  3.1902e-8\\
\hline
$(0.5,0.5)$  & $\rho_s$ & 0.9999  &  0.9926 &   0.6778 &  0.5173   & 0.5014  &  0.5004 &   0.4994  &  0.5000  &  0.5008
\\
             & $\rho_d$ & 4.2649e-1 &   8.0885e-3 &  1.3463e-4 &  8.6621e-6 &  8.1893e-7 &  8.1569e-8 &  8.1507e-9 &  8.1457e-10 &  8.1412e-11\\
\hline
$(0.9,0.9)$  & $\rho_s$ & 0.0849  &  0.4581  &  0.4950  &  0.5002  &  0.5003  &  0.5004  &  0.5005  &  0.5001 &   0.4998\\
             & $\rho_d$ & 3.1378e-1 &  1.8020e-1 &  2.1018e-2 &  2.1261e-3 &  2.1318e-4 &  2.1297e-5 &  2.1310e-6  & 2.1274e-7 &  2.1323e-8\\
\hline
$(0.99,0.99)$  & $\rho_s$ & 0.0003 &   0.0266  &  0.4140 &   0.4903 &   0.4986  &  0.4995  &  0.5005  &  0.4998   & 0.4999\\
            & $\rho_d$  & 2.5399e-3 &  2.4397e-3 &  2.1969e-3 &  3.1457e-4 &  3.2336e-5 &  3.2384e-6 &  3.2463e-7 &  3.2484e-8 &  3.2424e-9\\
\hline
\end{tabular}
\end{threeparttable}
\end{table*}
\endgroup

\begingroup
\renewcommand{\baselinestretch}{1}%
\begin{table*}[htbp!]
\centering
\begin{threeparttable}[b]
\renewcommand{\arraystretch}{1.3}
\caption{Statistical results of success rate and descent rate for the rotation transformation (Griewank problem)}
\scriptsize
\label{DSTA_benchmarks}
\begin{tabular}{{p{1.2cm}|p{0.6cm}|p{1.25cm}|p{1.25cm}|p{1.25cm}|p{1.25cm}|p{1.25cm}|p{1.25cm}|p{1.25cm}|p{1.25cm}|p{1.25cm}}}
\hline
$\mathrm{Best}_0$ & Index & $\alpha$ = 1 & $\alpha$ = 0.1 & $\alpha$ = 0.01 & $\alpha$ = 1e-3 & $\alpha$ = 1e-4 & $\alpha$ = 1e-5 & $\alpha$ = 1e-6 & $\alpha$ = 1e-7 & $\alpha$ = 1e-8\\
\hline
$(0.01,0.01)$ & $\rho_s$ &  0.0014 &   0.0967  &  0.4677  &  0.4966  &  0.4993 &   0.4996  &  0.5000  &  0.5007  &  0.5003\\
              & $\rho_d$ & 5.0031e-1 & 5.0462e-1 & 2.8800e-1 & 3.4337e-2 & 3.4921e-3 & 3.4973e-4 & 3.4962e-5  & 3.4957e-6 & 3.5030e-7\\
\hline
$(0.1,0.1)$  & $\rho_s$ & 0.0962 &   0.4676  &  0.4975  &  0.4997 &   0.5001  &  0.4992  &  0.5003  &  0.4995   & 0.4999\\
             & $\rho_d$ & 5.0587e-1 & 2.8807e-1 & 3.4234e-2 & 3.4888e-3 & 3.4825e-4&  3.4831e-5 & 3.4904e-6  & 3.4869e-7 & 3.4899e-8\\
\hline
$(0.5,0.5)$  & $\rho_s$ & 0.4250 &   0.4937 &   0.4993  &  0.4998  &  0.4996 &   0.5003  &  0.5000  &  0.5007  &  0.5006
\\
             & $\rho_d$ & 4.4869e-1 & 6.4092e-2&  6.5932e-3 & 6.6154e-4 & 6.6073e-5 & 6.6104e-6 & 6.6212e-7 &  6.6257e-8 & 6.6127e-9\\
\hline
$(0.9,0.9)$  & $\rho_s$ & 0.4532 &   0.4958 &   0.4994 &   0.5005 &   0.4998  &  0.4995 &   0.4996  &  0.5002  &  0.4994\\
             & $\rho_d$ & 2.7822e-1 &  3.1695e-2 & 3.2066e-3 & 3.2065e-4 & 3.2064e-5 & 3.2016e-6 & 3.2054e-7  & 3.2012e-8 & 3.2048e-9\\
\hline
$(0.99,0.99)$  & $\rho_s$ & 0.4550 &   0.4958  &  0.4999  &  0.5003  &  0.5002 &   0.4996 &   0.4999  &  0.5003  &  0.5002\\
            & $\rho_d$  & 2.4931e-1 & 2.7546e-2 & 2.7818e-3 & 2.7817e-4 & 2.7754e-5 & 2.7808e-6 & 2.7839e-7  & 2.7824e-8 & 2.7801e-9\\
\hline
\end{tabular}
\end{threeparttable}
\end{table*}
\endgroup

\begin{figure*}[!htbp]
  \includegraphics[width=6.5cm,height=5cm]{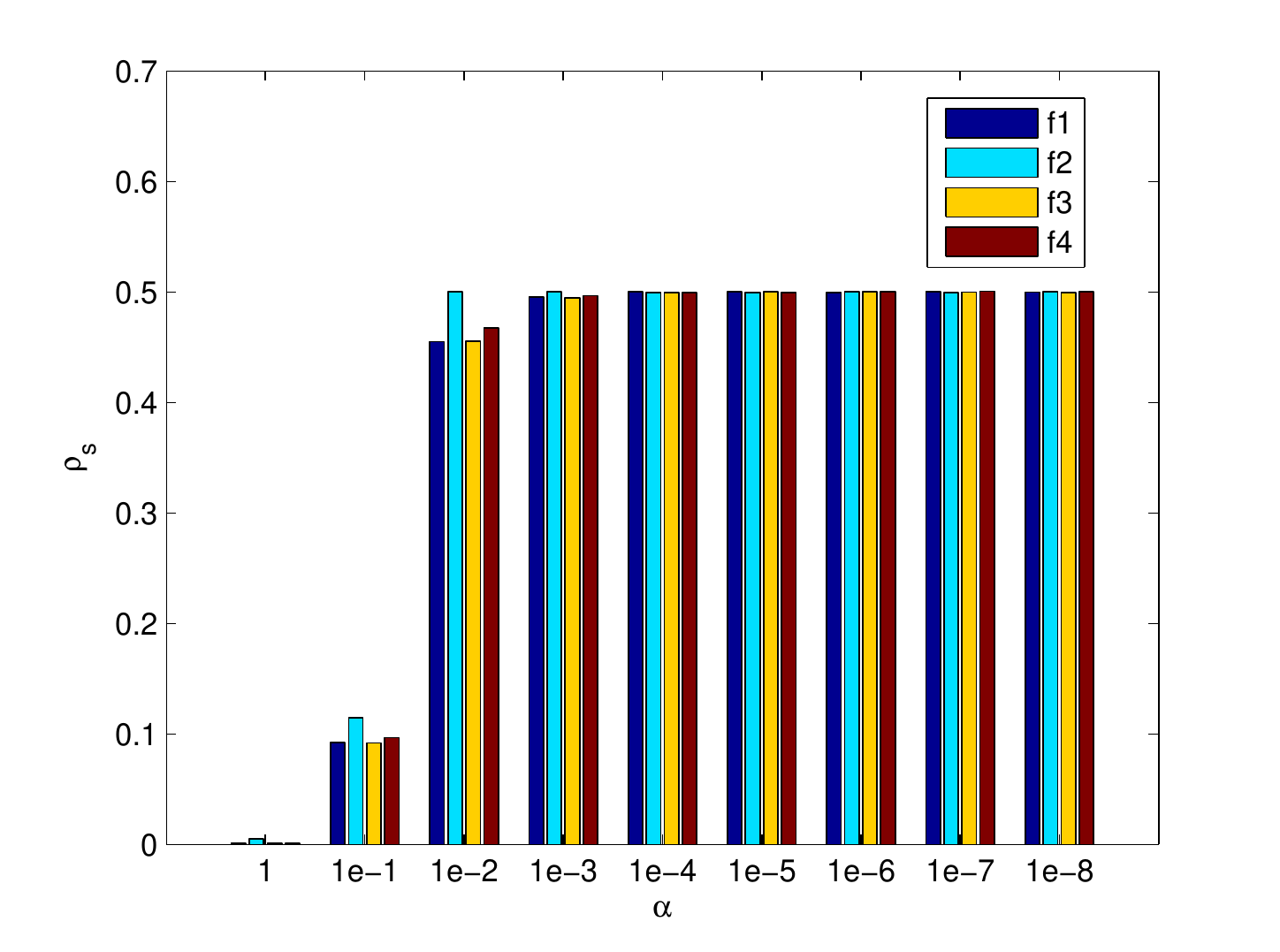}
  \includegraphics[width=6.5cm,height=5cm]{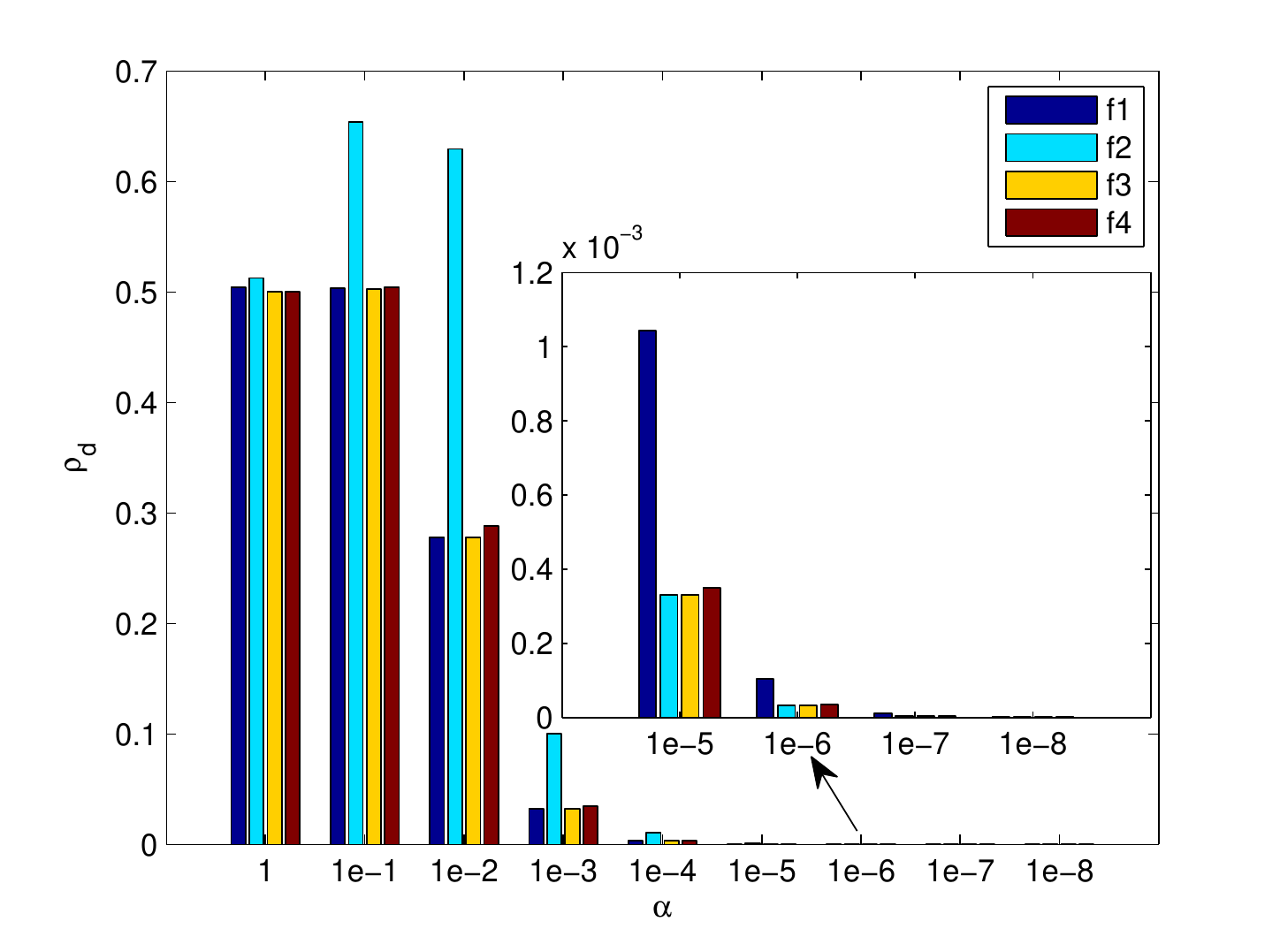}
  \includegraphics[width=6.5cm,height=5cm]{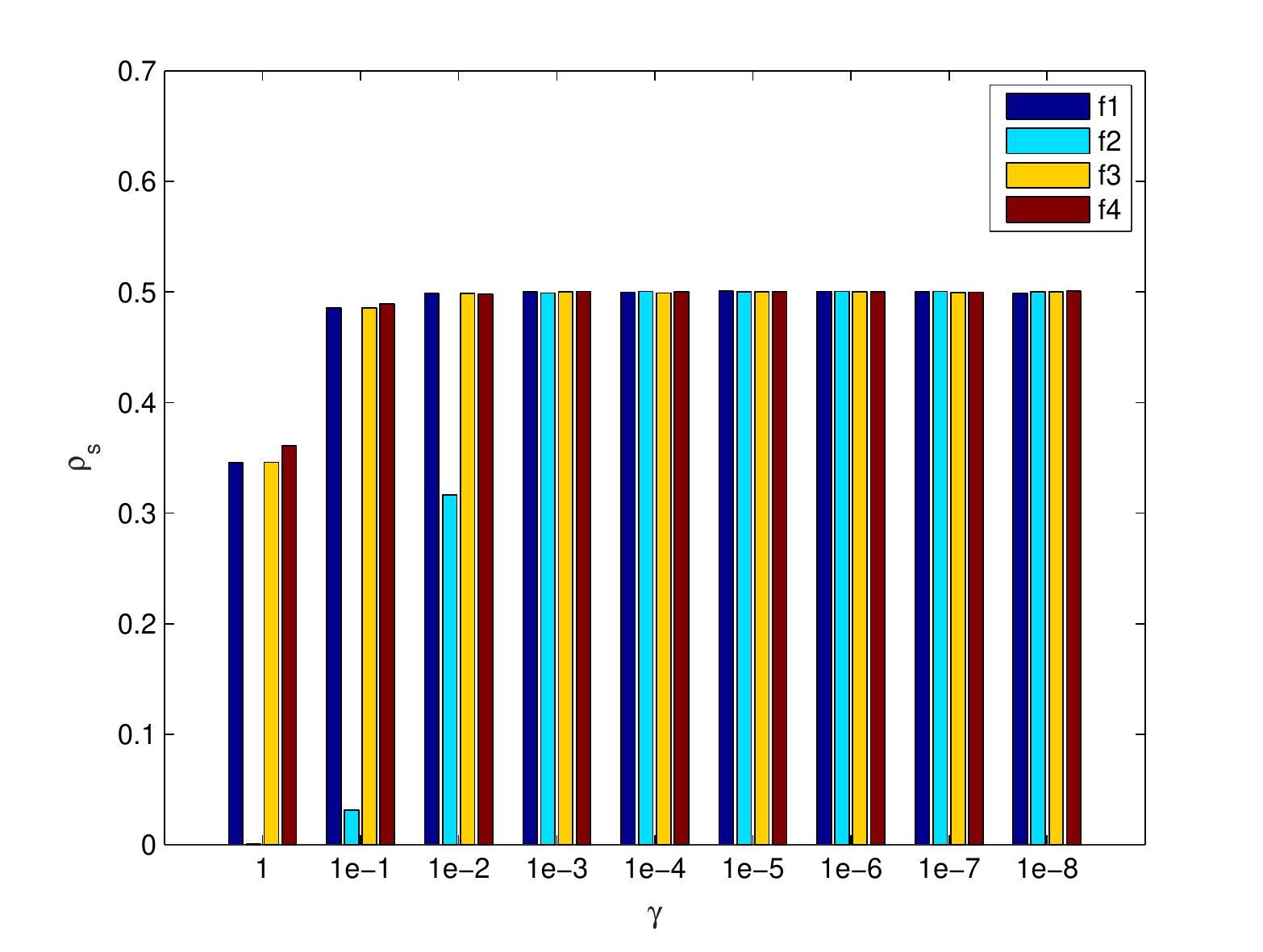}
  \includegraphics[width=6.5cm,height=5cm]{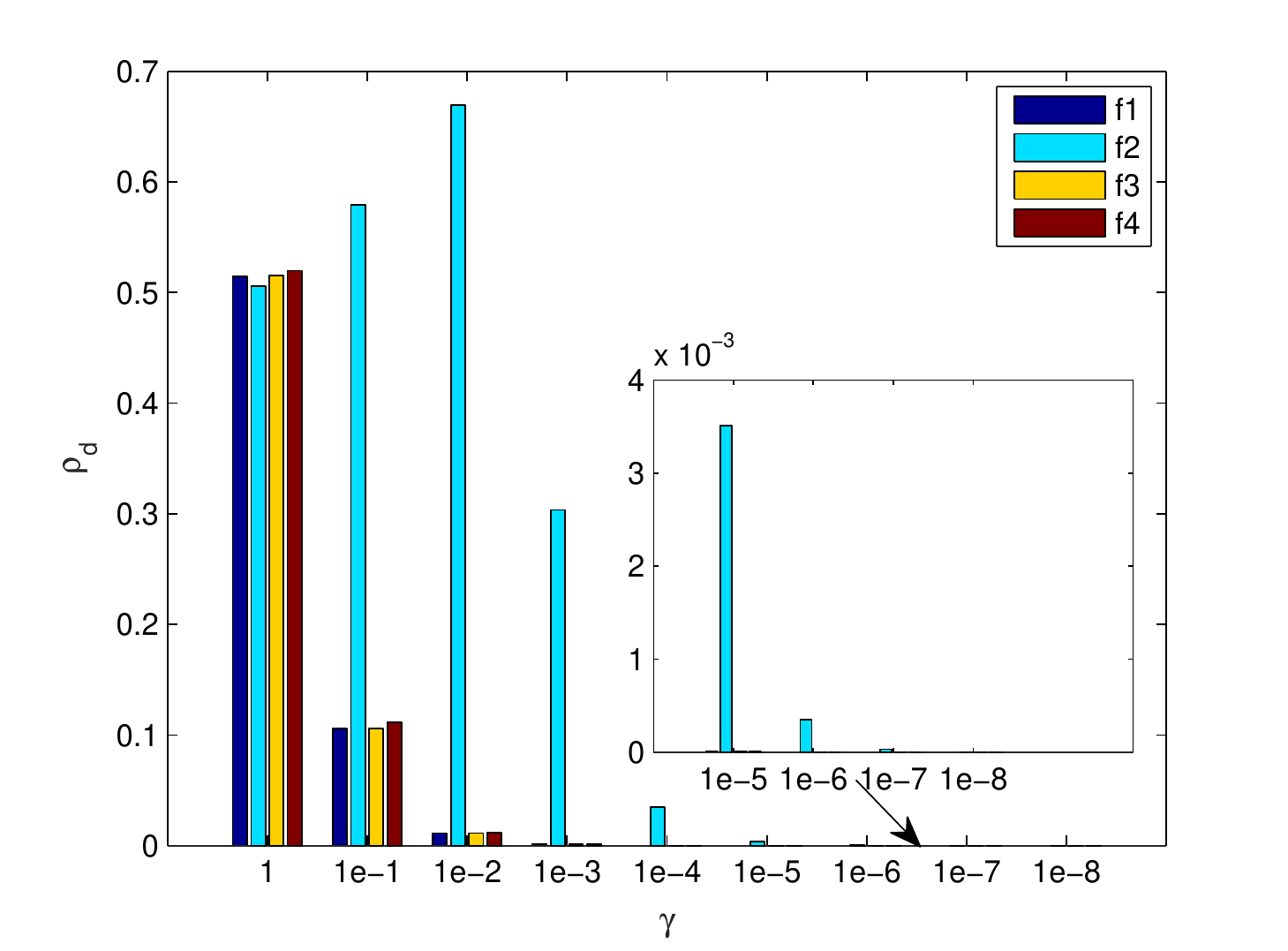}
  \includegraphics[width=6.5cm,height=5cm]{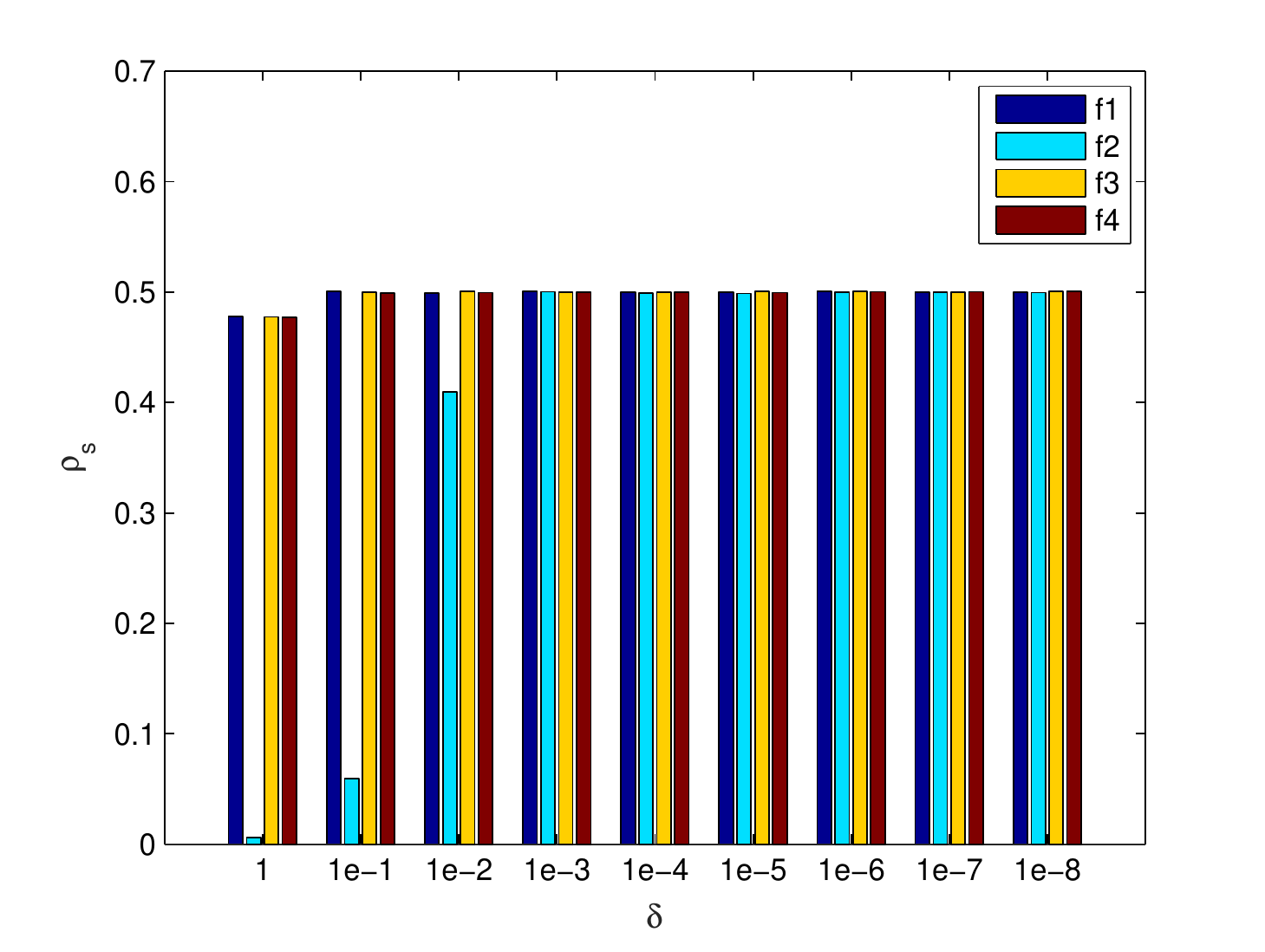}
  \includegraphics[width=6.5cm,height=5cm]{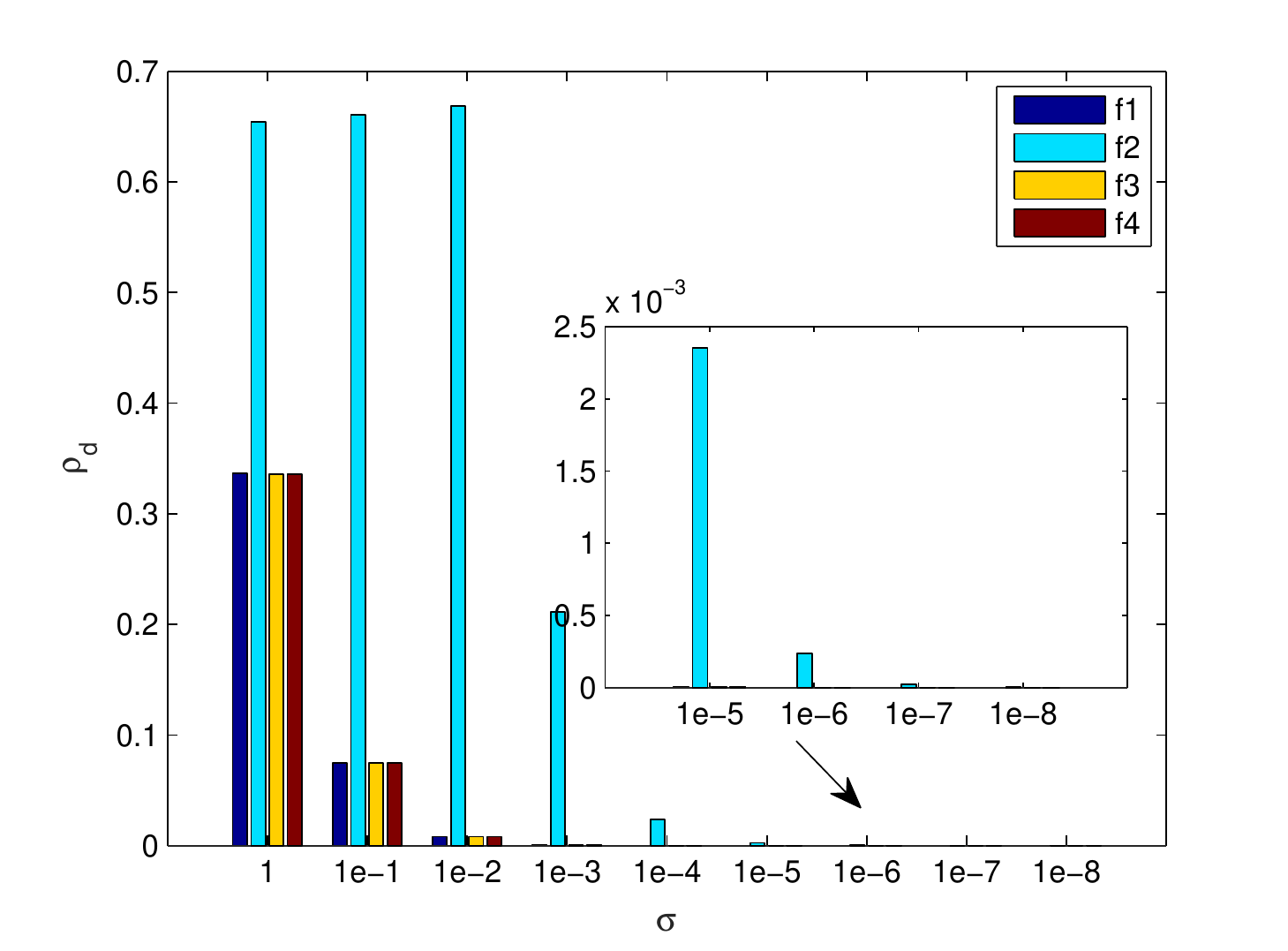}
\caption{The changes of success rate $\rho_s$ and descent rate $\rho_d$ with
 $\alpha$, $\gamma$ and $\delta$ respectively when approaching the global minima}
\label{fig_expansion}       
\end{figure*}

As indicated in these tables, the following properties can be observed:
\begin{enumerate}
  \item as the decrease of a state transformation factor below a certain threshold, the descent rate $\rho_d$ is showing
  a declining trend.
  \item the success rate remains almost steadily high if a state transformation factor is below a certain threshold.
  \item the success rate of the rotation transformation is not high until the rotation factor is below a
  threshold when current solution $Best_0$
  is approaching the global optimal solution.
\end{enumerate}

To be more specific, let's take the rotation transformation for example, the changes of success rate $\rho_s$ and
descent rate $\rho_d$ with the rotation factor $\alpha$ are illustrated in Fig. 1 when current solution
$Best_0$
is approaching the global optimum. Here, $Best_0$ equals to (0.01, 0.01), (0.99, 0.99), (0.01, 0.01) and (0.01, 0.01)
for $f_1$, $f_2$, $f_3$ and $f_4$ respectively.
By taking a closer look at these two figures, it is not difficult
to find that there exists a trade-off between the success rate and the descent rate.
For instance, when $\alpha = 1 $, the success rate $\rho_s$ is quite low, while the descent rate $\rho_d$ is quite high.
On the contrary, when $\alpha \in \{$1e-5, 1e-6, 1e-7, 1e-8$\}$, the success rate $\rho_s$ is quite high,
while the descent rate $\rho_d$ is quite low.
\begin{remark}
The property 3) can provide additional support to the way in changing the rotation factor in the standard
continuous STA, \textit{i.e.}, $\alpha$ is not kept constant but decreasing periodically from a maximum
value $\alpha_{\max}$ to a minimum value $\alpha_{\min}$. Anyway, it is obvious that the way in changing the state transition factors is not in an optimal manner.
\end{remark}

\section{State transition algorithm with optimal parameter selection}
As inspired by the statistical study of the state transformation factors, in this section,
an optimal parameter selection strategy is proposed
to accelerate the search of the standard continuous state transition algorithm.
\subsection{Optimal parameter selection for the state transformation factors}
In classical iterative methods for numerical optimization, the following iterative formula is usually adopted
\begin{eqnarray}\label{eq_numopt}
\bm x_{k+1} = \bm x_k + a_k \bm d_k
\end{eqnarray}
where $\bm d_k$ is the search direction, and $a_k$ is the step size.
For gradient-based algorithms, the search direction is relevant to the gradient of current iterative point, for instance, the steepest descent method, $\bm d_k = -\nabla f(\bm x_k)$,
and the step size $a_k$ is often restricted to the range [0,1].
It can be found that the pattern of iterative formula in continuous STA is similar
to that of Eq. (\ref{eq_numopt}), as shown below
\begin{displaymath}
\left. \begin{array}{l}
\frac{1}{n \|\bm s_{k}\|_{2}} R_{r} \bm s_{k} \\
R_{t} \frac{\bm s_{k}- \bm s_{k-1}}{\|\bm s_{k}- \bm s_{k-1}\|_{2}}\\
R_{e} \bm s_k \\
R_{a} \bm s_k
\end{array} \right \} \Rightarrow \tilde{\bm d}_k, \;\;\;
\left. \begin{array}{l}
\alpha \\
\beta \\
\gamma \\
\delta
\end{array} \right \} \Rightarrow \tilde{a}_k,
\end{displaymath}
and a big difference is that
the search direction is not determined.
Compared with gradient-based algorithms, the STA can be used for global optimization lies in at least two aspects: 1)
the search is in all directions; 2) the search can go to any length. While compared with the traditional trust region method, the similarity is that
some parts of the STA (except the translation transformation) can be considered as a special kind of trust region method, but the differences are: i) the STA utilizes the original function not its quadratic approximation; ii) the search direction in STA is stochastic.

For rotation and translation transformation, the search zone is restricted in a hypersphere or along a line, which are controlled by the corresponding transformation factors. For expansion and axesion transformation, although the search zone can be expanded to the whole space in probability due to the Gaussian distribution, the search zone is restricted and manipulated by the expansion and axesion factors as well. That is to say, in practical numerical computation,
the neighborhood formed by the state transformation operators is controlled by the transformation factors to a large extent, which are also testified by the statistical study. To simplify the parameter selection and accelerate the search process, the values of these parameters are all taken from the
set $\Omega = \{$1, 1e-1, 1e-2, 1e-3, 1e-4, 1e-5, 1e-6, 1e-7, 1e-8$\}$, and the parameter value
with the corresponding smallest objective function value is chosen.

\subsection{The proposed STA}
Let's denote the optimal parameter as $\tilde{a}^{*}$, and then we have
\begin{eqnarray}
\tilde{a}^{*} = \argmin_{\tilde{a}_k \in \Omega} f(\bm x_k + \tilde{a}_k \tilde{\bm d_k})
\end{eqnarray}

In theory, the neighborhood formed by the state transformation operators has infinite candidate solutions; however, only \textit{SE} samples are
used for evaluation in practice. That is to say, for a given parameter value, only \textit{SE} samples are taken into consideration.
In order to further utilize the parameter more completely, the selected parameter value is kept for a period of time, denoted as $T_p$.
To be more specific, the detailed of the proposed STA can be outlined as follows
\begin{algorithmic}[1]
\Repeat
    \State {Best $\gets$ expansion\_w(funfcn,Best,SE,$\Omega$)}
    \State {Best $\gets$ rotation\_w(funfcn,Best,SE,$\Omega$)}
    \State {Best $\gets$ axesion\_w(funfcn,Best,SE,$\Omega$)}
\Until{the specified termination criterion is met}
\end{algorithmic}

\indent In the meanwhile, rotation\_w$(\cdot)$ in above pseudocode is given for further explanations
\begin{algorithmic}[1]
\State{[Best,$\alpha$] $\gets$ update\_alpha(funfcn,Best,SE,$\Omega$)}
\For{$i\gets 1, T_p$}
\State {Best $\gets$ rotation(funfcn,Best,SE,$\alpha$)}
\EndFor
\end{algorithmic}
where the function update\_alpha represents the implementation of selection the optimal parameter value of rotation factor. The proposed STA differs from the standard STA in three folds: 1) the periodical way of diminishing the transformation factors is no longer used;
2) the optimal parameter is selected for state transformation; 3) the optimal parameter is kept to utilize for a period of time.

\section{Experimental results}
In order to testify the effectiveness of the proposed STA, the following additional benchmark functions are used for test.

(5) Ackley function
\[
f_5(\bm x)\!=\!20\!+\!e\!-\!20\exp(\!-\!0.2\sqrt{\frac{1}{n} \sum_{i=1}^n x_{i}^2})
\!-\!\exp(\frac{1}{n} \sum_{i=1}^n \cos(2\pi x_{i}))
\]
where the global optimum $\bm x^{*} = (0, \cdots, 0)$ and $f(\bm x^{*}) = 0$,
$-32 \leq x_i \leq 32, i= 1, \cdots, n$.

(6) High Conditioned Elliptic function
\[
f_6(\bm x) = \sum_{i=1}^n (10^6)^{\frac{i-1}{n-1}}x_i^2
\]
where the global optimum $\bm x^{*} = (0, \cdots, 0)$ and $f(\bm x^{*}) = 0$,
$-100 \leq x_i \leq 100, i= 1, \cdots, n$.

(7) Michalewicz function
\[
f_7(\bm x) = -\sum_{i=1}^n \sin(x_{i}) \sin( \frac{i x_{i}^2}{\pi})^{20}
\]
where the global optimum is unknown,
$0 \leq x_i \leq \pi, i= 1, \cdots, n$.

(8) Trid function
\[
f_8(\bm x) = \sum_{i=1}^n (x_i - 1)^2 -\sum_{i=2}^n x_i x_{i-1}
\]
where the global optimum $x_i^{*} = i(n+1-i)$ and $f(\bm x^{*}) = -\frac{n(n+4)(n-1)}{6}$,
$-n^2 \leq x_i \leq n^2, i= 1, \cdots, n$.

(9) Schwefel function
\[
f_9(\bm x) = \sum_{i=1}^n [ -x_{i} \sin(\sqrt{|x_{i}|})]
\]
where the global optimum $\bm x^{*} = (420.9687, \cdots, 420.9687)$ and $f(\bm x^{*}) = -418.9829n$,
$-500 \leq x_i \leq 500, i= 1, \cdots, n$.

(10) Schwefel 1.2 function
\[
f_{10}(\bm x) = \sum_{i=1}^n (\sum_{j=1}^i x_j)^2
\]
where the global optimum $\bm x^{*} = (0, \cdots, 0)$ and $f(\bm x^{*}) = 0$,
$-100 \leq x_i \leq 100, i= 1, \cdots, n$.

(11) Schwefel 2.4 function
\[
f_{11}(\bm x) = \sum_{i=1}^n [(x_i - 1)^2 + (x_1 - x_i^2)^2
\]
where the global optimum $\bm x^{*} = (1, \cdots, 1)$ and $f(\bm x^{*}) = 0$,
$0 \leq x_i \leq 10, i= 1, \cdots, n$.

(12) Weierstrass function
\[
f_{12}(\bm x) = \sum_{i=1}^n\sum_{k=0}^{k_{\max}}[a^k\cos(2 \pi b^k(x_i+0.5))] - n \sum_{k=0}^{k_{\max}}a^k\cos(\pi b^k x_i)
\]
where $a = 0.5, b = 3, k_{\max} = 20$, the global optimum $\bm x^{*} = (0, \cdots, 0)$ and $f(\bm x^{*}) = 0$,
$-0.5 \leq x_i \leq 0.5, i= 1, \cdots, n$.

\begin{table*}[!htbp]
\begin{center}
\caption{Comparisons among various algorithms on test functions}
\scriptsize
\begin{tabular}{{p{0.1cm}p{0.1cm}p{2.65cm}p{2.65cm}p{2.6cm}p{2.6cm}p{2.6cm}p{2.85cm}}}
\hline
\toprule[1pt]
Fcn & Dim & GL-25 & CLPSO & SaDE & ABC & Standard STA & Proposed STA \\
\hline
               & 20 & 2.5523e-10 $\pm$ 1.5883e-10  & 6.2546e-43 $\pm$ 1.4407e-42 &6.3533e-188 $\pm$ 0 &2.7287e-16 $\pm$ 6.1809e-17 & 0 $\pm$ 0 & 0 $\pm$ 0	 \\
$f_{1} $       & 30 & 1.7872e-8  $\pm$ 1.0381e-8  &1.9944e-40 $\pm$ 1.8175e-40  &5.5498e-184 $\pm$ 0&5.6618e-16 $\pm$ 7.3169e-17& 0 $\pm$ 0 & 0 $\pm$ 0	   \\
               & 50 & 2.3336e-6  $\pm$ 1.5613e-6   &8.3697e-63 $\pm$ 6.6314e-63  &4.8082e-190 $\pm$ 0 &1.3115e-15 $\pm$ 1.4686e-16& 0 $\pm$ 0  & 0 $\pm$ 0  \\
\hline
               & 20 & 15.9120 $\pm$ 0.2273 &1.3524 $\pm$ 1.5792&0.7973 $\pm$ 1.6361 &0.0871 $\pm$ 0.1254 & 0.0327 $\pm$ 0.0019 & 3.2981e-07 $\pm$	 1.0312e-06  \\
$f_{2}$        & 30 & 25.9785  $\pm$ 0.1774   & 3.3395 $\pm$ 4.4690& 1.3895 $\pm$ 2.1499&0.0523 $\pm$ 0.0672& 0.0711 $\pm$ 0.0128 & 1.0027e-07 $\pm$ 1.0502e-07   \\
               & 50 & 46.3067  $\pm$ 0.4004   & 38.4515 $\pm$ 31.7815 & 16.2265 $\pm$ 21.2962&0.0634 $\pm$ 0.1142& 2.5228 $\pm$ 1.2541 & 1.0660e-07 $\pm$ 8.0190e-08  \\
\hline
               & 20 & 88.3377 $\pm$ 10.1747   & 0 $\pm$ 0 & 0.2985 $\pm$ 0.4678&4.2633e-15 $\pm$ 1.0412e-14& 0 $\pm$ 0	 & 0 $\pm$ 0 \\
$f_{3} $       & 30 & 177.1109  $\pm$ 12.2431   &5.6843e-15 $\pm$ 1.7496e-14 &1.0945 $\pm$ 0.8479 & 8.5265e-14 $\pm$ 4.3251e-14& 0 $\pm$ 0  & 0 $\pm$ 0  \\
               & 50 & 365.4491  $\pm$ 12.8696   &0 $\pm$ 0 & 5.5220 $\pm$ 2.6516 &1.0601e-12 $\pm$ 1.6704e-12 & 0 $\pm$ 0 &8.5265e-15 $\pm$ 2.7817e-14        \\
\hline
               & 20 &  0.2620 $\pm$ 0.1020   &  5.2736e-16 $\pm$ 1.4066e-15 &0.0034  $\pm$ 0.0051 &1.3711e-15 $\pm$ 2.2887e-15&0 $\pm$ 0 &0 $\pm$ 0 \\
$f_{4}$        & 30 & 0.0178  $\pm$ 0.0797   &0 $\pm$ 0  & 0.0041 $\pm$ 0.0098 &7.9936e-16 $\pm$ 6.3277e-16 & 0 $\pm$ 0 &0 $\pm$ 0\\
               & 50 & 2.0621e-6  $\pm$ 1.2609e-6   &0 $\pm$ 0  & 0.0229 $\pm$ 0.0338 &1.5432e-15 $\pm$ 6.5721e-16 & 0 $\pm$ 0 &0 $\pm$ 0  \\
\hline
               & 20 &  2.9519e-6 $\pm$ 8.5896e-7   &6.0396e-15 $\pm$ 7.9441e-16 & 2.6645e-15 $\pm$ 0&2.4336e-14 $\pm$ 3.6267e-15&7.1054e-16 $\pm$  1.8134e-15& 1.2434e-15 $\pm$ 1.7857e-15  \\
$f_{5} $       & 30 & 1.7312e-5  $\pm$ 4.9711e-6   &7.2831e-15 $\pm$ 1.6704e-15& 0.4004 $\pm$ 0.5677 &4.7073e-14 $\pm$ 5.2189e-15&  2.4869e-15 $\pm$ 7.9441e-16 & 2.6645e-15 $\pm$ 0 \\
               & 50 & 1.5638e-4  $\pm$ 5.1830e-5   &1.2790e-14 $\pm$ 1.3015e-15& 1.8811 $\pm$ 1.8811 & 1.0214e-13 $\pm$ 8.3914e-15&2.6645e-15 $\pm$ 0 &2.6645e-15 $\pm$ 0 \\
\hline
               & 20 & 1.0920e-7 $\pm$ 8.5896e-7   &9.6865e-40 $\pm$ 1.0089e-39&2.9962e-182 $\pm$ 0  &2.8100e-16 $\pm$ 2.2693e-17 & 0 $\pm$ 0 & 0 $\pm$ 0	  \\
$f_{6}$        & 30 & 4.3319e-6  $\pm$ 4.1372e-6   &9.6865e-40 $\pm$ 1.0089e-39& 7.5626e-180 $\pm$ 0 &5.0197e-16 $\pm$ 5.0710e-17& 0 $\pm$ 0  & 0 $\pm$ 0 \\
               & 50 & 1.4938e-4  $\pm$ 9.1548e-5   & 1.0767e-59 $\pm$ 9.2087e-60& 3.6584e-188 $\pm$ 0 &1.2513e-15 $\pm$ 1.3320e-16&0 $\pm$ 0  & 0 $\pm$ 0\\
\hline
               & 20 &  -10.7121 $\pm$ 0.4311   & -19.6363 $\pm$ 0.0013 & -19.6204 $\pm$ 0.0210  &-19.6359 $\pm$ 0.0013&-19.2512 $\pm$ 0.7144& -19.6370 $\pm$ 4.5865e-15 \\
$f_{7} $       & 30& -13.5080 $\pm$ 4.1372e-6 & -29.5405 $\pm$ 0.0422 & -29.5668 $\pm$ 0.0439 &-29.6083 $\pm$ 0.0121& -29.2917 $\pm$ 0.5761 & -29.3322 $\pm$ 0.4810 \\
               & 50 & -18.2114  $\pm$ 0.8100   &   -49.2281 $\pm$ 0.1068&-49.3694 $\pm$ 0.1155 &-49.5258 $\pm$ 0.0239&-48.9364 $\pm$ 0.7706&  -49.2284 $\pm$ 0.5118\\
\hline
               & 20 & -1.2099e3 $\pm$  63.9563   &  -1.2126e3 $\pm$ 198.0180  & -1.5200e3 $\pm$ 0.0030 &-1.4934e3 $\pm$ 17.7408&-1.5200e3 $\pm$  7.6568e-10& -1.5200e3 $\pm$ 1.0526e-09  \\
$f_{8} $       & 30 & -2.1886e3 $\pm$  419.3236   &-2.3303e3 $\pm$ 786.8715 & -4.8684e3 $\pm$ 30.0984  &-3.7616e3 $\pm$ 379.9890& -4.9300e3 $\pm$ 1.0317e-8& -4.9300e3 $\pm$ 1.9630e-8                \\
               & 50 & -3.8943e3  $\pm$ 1.0972e3   & -4.6039e3 $\pm$ 2.1276e3& -1.6988e4 $\pm$ 1.9831e3 &-6.2731e3 $\pm$ 2.9887e3 &-2.2050e4 $\pm$ 2.4728e-5 & -2.2050e4 $\pm$ 1.2703e-06       \\
\hline

               & 20 & -3.4543e3 $\pm$ 262.3109   & -8.3797e3 $\pm$ 3.7325e-12&  -8.3678e3 $\pm$ 52.9672 & -8.3797e3 $\pm$ 1.7705e-12 & -8.3797e3 $\pm$ 2.0013e-12 &-8.3797e3 $\pm$ 2.4688e-12	                   \\
$f_{9}$        &30 &  -4.2340e3  $\pm$ 206.4148   &  -1.2569e4 $\pm$ 1.8662e-12& -1.2534e4 $\pm$ 55.6852&-1.2569e4 $\pm$ 5.0878e-10 &  -1.2569e4 $\pm$  3.9808e-12 & -1.2569e4 $\pm$ 3.9808e-12	                 \\
               & 50 &-5.5094e3  $\pm$ 308.3849   &-2.0949e4 $\pm$ 7.4650e-12&-2.0848e4 $\pm$ 123.1747&-2.0949e4 $\pm$ 2.2694e-4&-2.0949e4 $\pm$ 6.9328e-12 &-2.0949e4 $\pm$ 6.8824e-12        \\
\hline
               & 20 & 680.7399 $\pm$  158.9219   &6.0736 $\pm$ 3.4399 &3.4609e-29 $\pm$ 1.5188e-28 &122.8894 $\pm$ 57.5833& 0 $\pm$ 0 & 5.4347e-323 $\pm$ 0	                   \\
$f_{10}$       & 30 &  6.0084e3  $\pm$ 891.0211   & 192.7349 $\pm$ 40.2446 &1.7348e-18 $\pm$ 3.7702e-18&1.4779e3 $\pm$ 417.2093
& 0 $\pm$ 0&2.9857e-207 $\pm$ 0                 \\
               & 50 & 2.6141e4  $\pm$ 3.4905e3   &1.7640e3 $\pm$ 334.6602& 1.1549e-11 $\pm$ 2.5786e-11&1.3111e4 $\pm$ 2.0672e3& 4.9187e-16 $\pm$ 2.1481e-15& 1.3290e-143 $\pm$ 2.4407e-143    \\
\hline
               & 20 & 2.6106e-9 $\pm$ 2.0321e-9   &5.9122e-6 $\pm$ 2.7178e-6& 1.8933e-30 $\pm$ 1.2854e-30 &4.6418e-15 $\pm$ 1.2107e-15 & 3.4527e-13  $\pm$  1.9644e-13 &2.8663e-21 $\pm$ 1.5596e-21	                   \\
$f_{11}$       & 30 & 1.8994e-9 $\pm$ 1.5212e-9  &7.2339e-5 $\pm$ 1.6458e-5&6.9050e-30 $\pm$ 1.9571e-30& 7.9455e-15 $\pm$ 1.9605e-15&4.3734e-13  $\pm$  2.3603e-13& 4.7017e-21 $\pm$ 1.8546e-21                 \\
               & 50 & 6.2123e-12  $\pm$ 1.0570e-11   &4.6982e-7 $\pm$ 1.0856e-7& 4.1750e-29 $\pm$ 2.8959e-29&2.3607e-14 $\pm$ 7.8298e-15&8.2000e-13 $\pm$ 2.3292e-13 & 8.3954e-21 $\pm$ 2.8656e-21       \\
\hline
               & 20 & 0.0047 $\pm$ 0.0012   & 0 $\pm$ 0 & 0 $\pm$ 0 &0 $\pm$ 0& 0 $\pm$ 0 & 0 $\pm$ 0                  \\
$f_{12}$       & 30 &0.0302 $\pm$ 0.0077 & 0 $\pm$ 0 & 0.2097 $\pm$ 0.3033 &0 $\pm$ 0& 0 $\pm$ 0 & 0 $\pm$ 0	                 \\
               & 50 & 0.2307  $\pm$  0.0805   & 0 $\pm$ 0 & 1.7162 $\pm$ 0.8589 &2.7711e-14 $\pm$ 9.7534e-15& 0 $\pm$ 0  & 0 $\pm$ 0      \\
\bottomrule[1pt]
\hline
\end{tabular}
\end{center}
\end{table*}

\begin{figure*}[!htbp]
  \includegraphics[width=9.3cm,height=7cm]{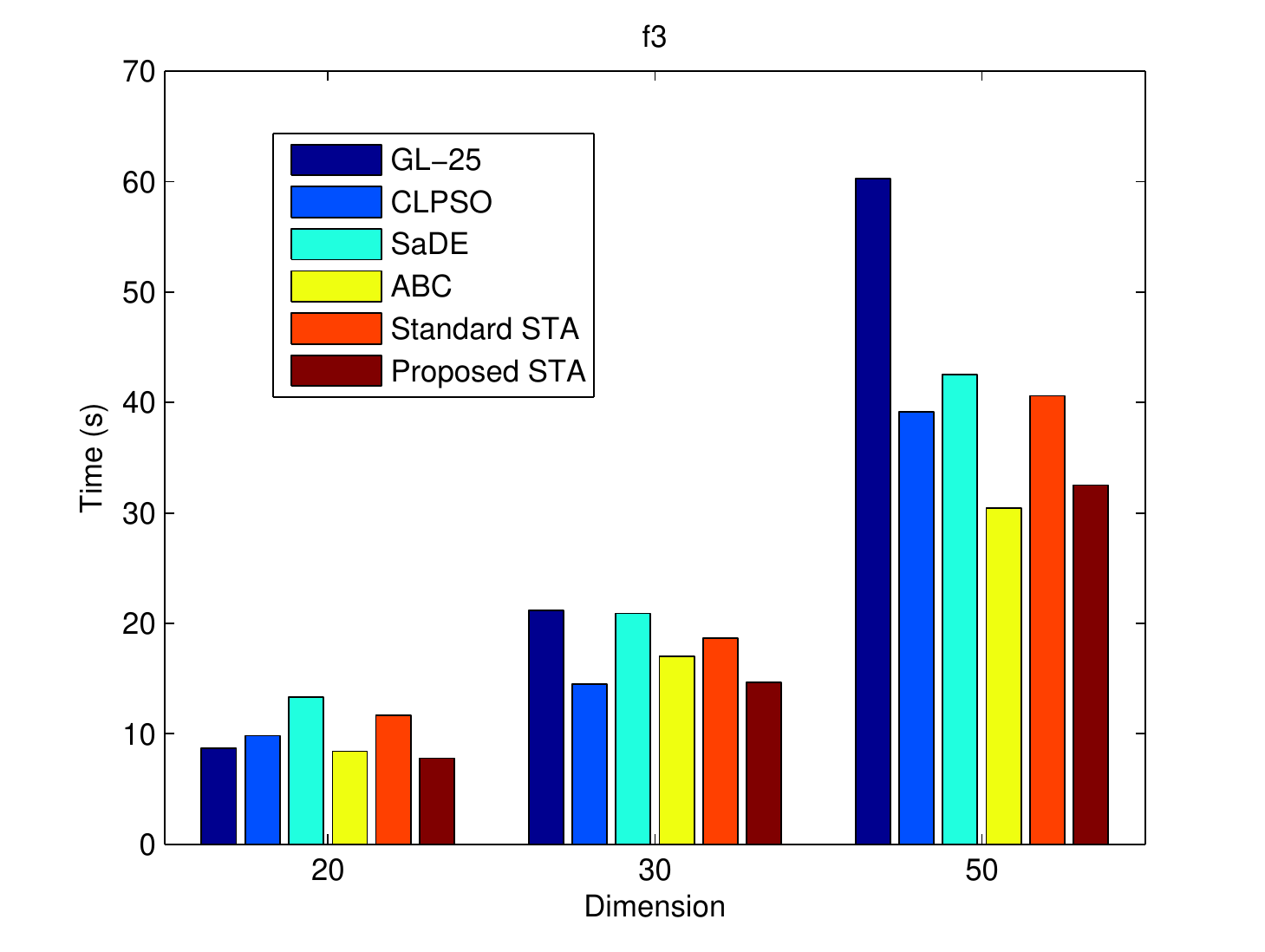}
  \includegraphics[width=9.3cm,height=7cm]{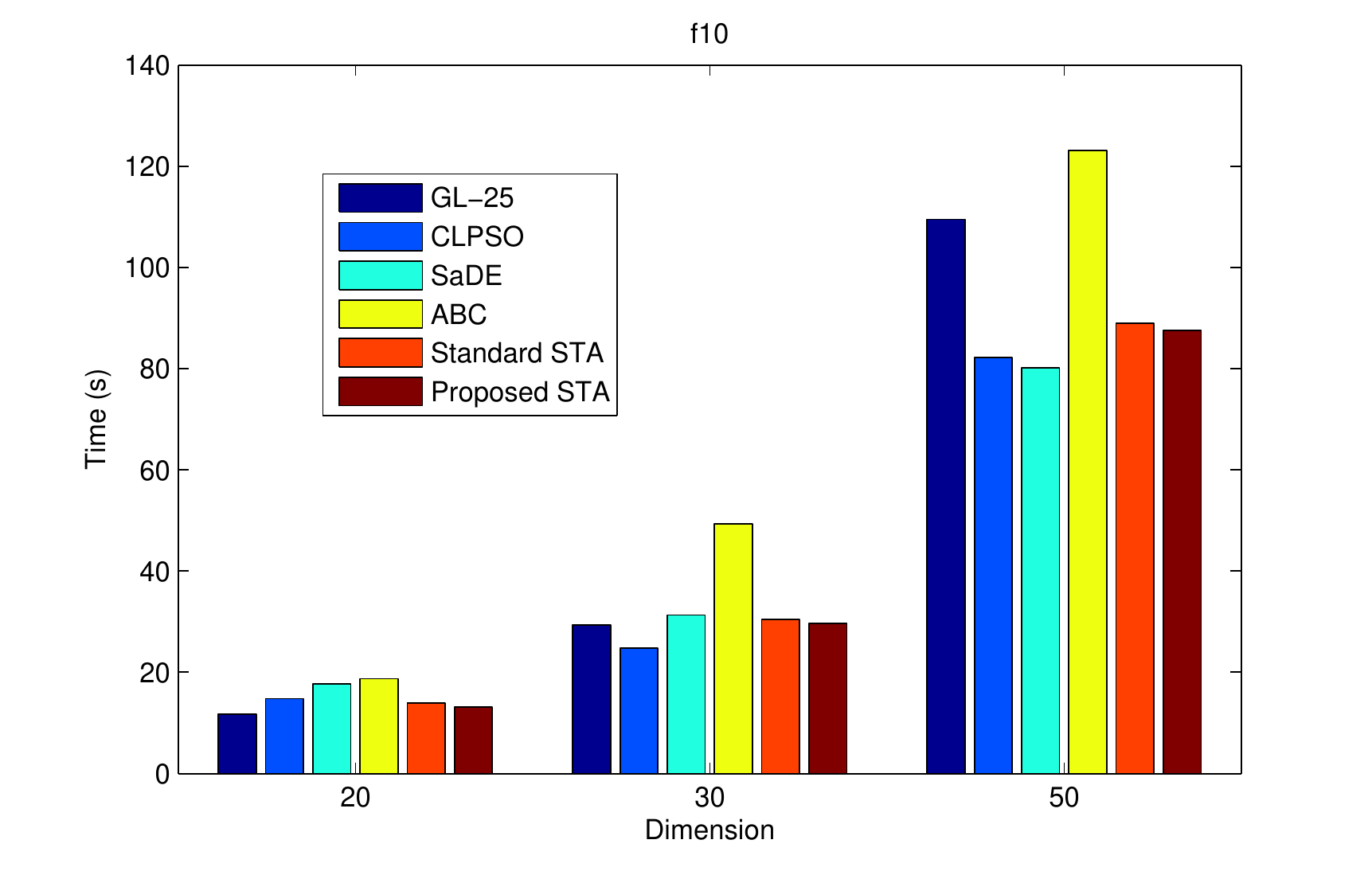}
\caption{The average elapsed time for different metaheuristic methods with respect to $f_3$ and $f_{10}$, respectively}
\label{fig_expansion}       
\end{figure*}

\begin{figure*}[!htbp]
  \includegraphics[width=6.5cm,height=5cm]{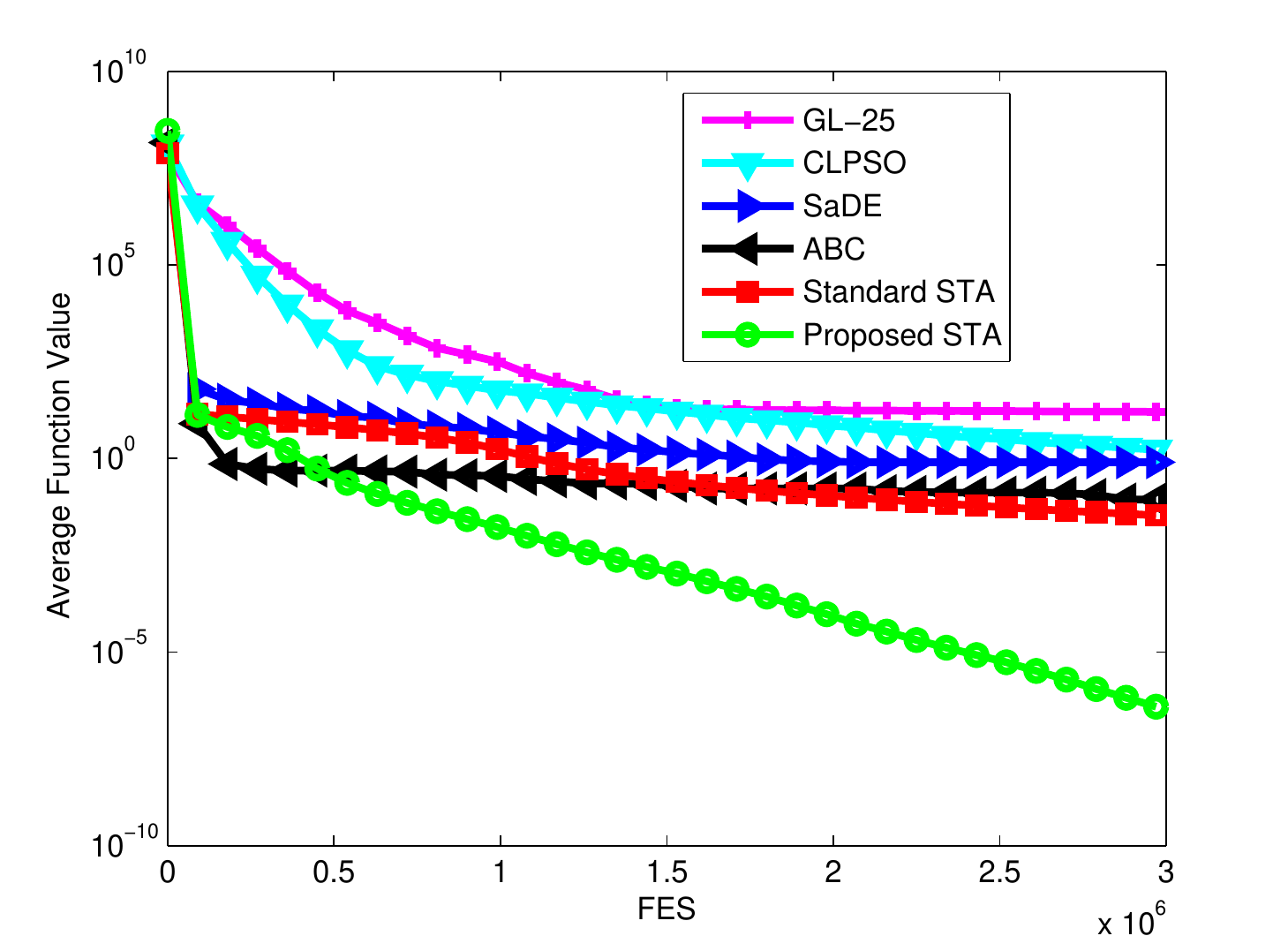}
  \includegraphics[width=6.5cm,height=5cm]{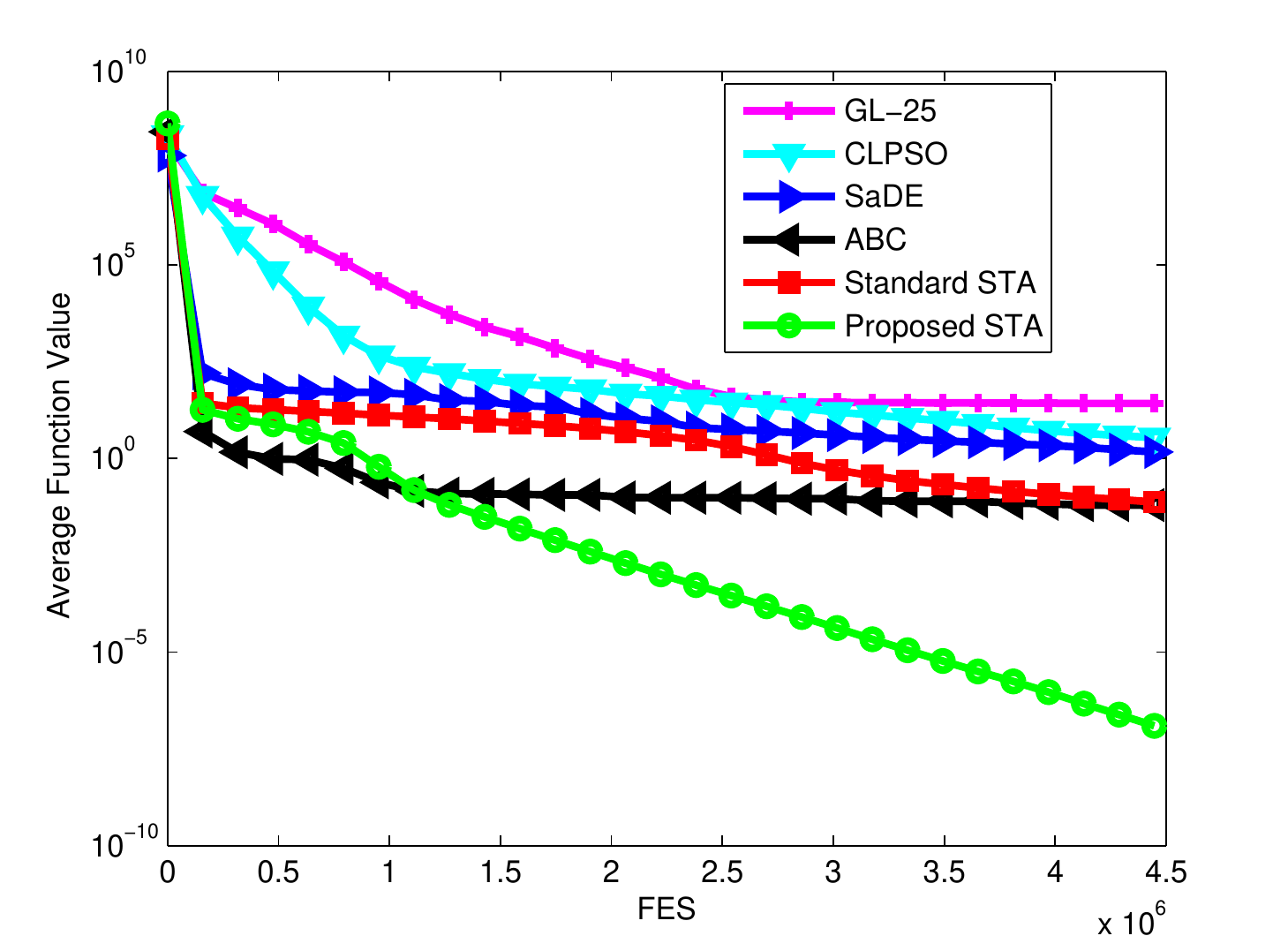}
  \includegraphics[width=6.5cm,height=5cm]{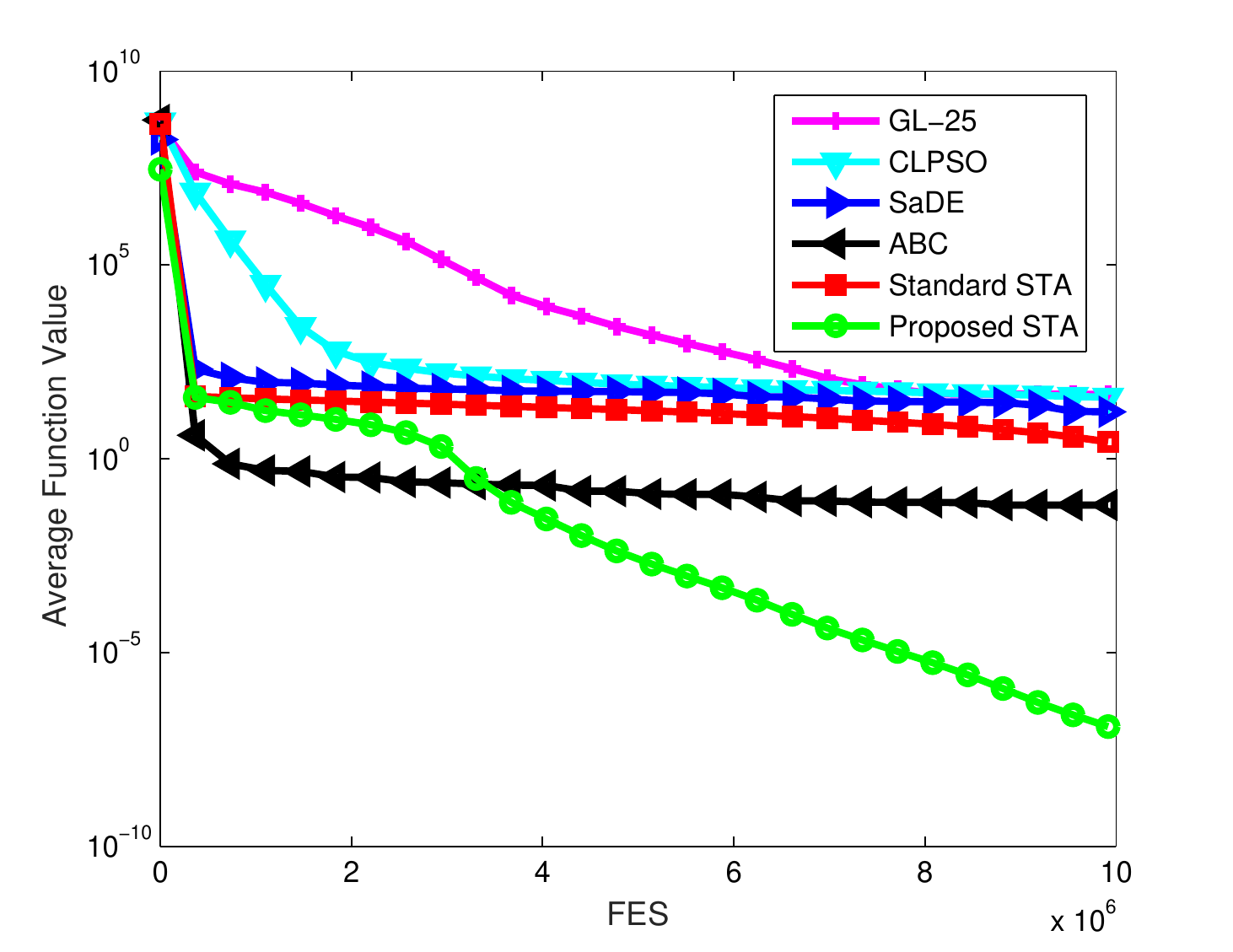}
\caption{The average iterative curves for different metaheuristic methods with respect to the Rosenbrock function}
\label{fig_expansion}       
\end{figure*}

\begin{figure*}[!htbp]
  \includegraphics[width=6.5cm,height=5cm]{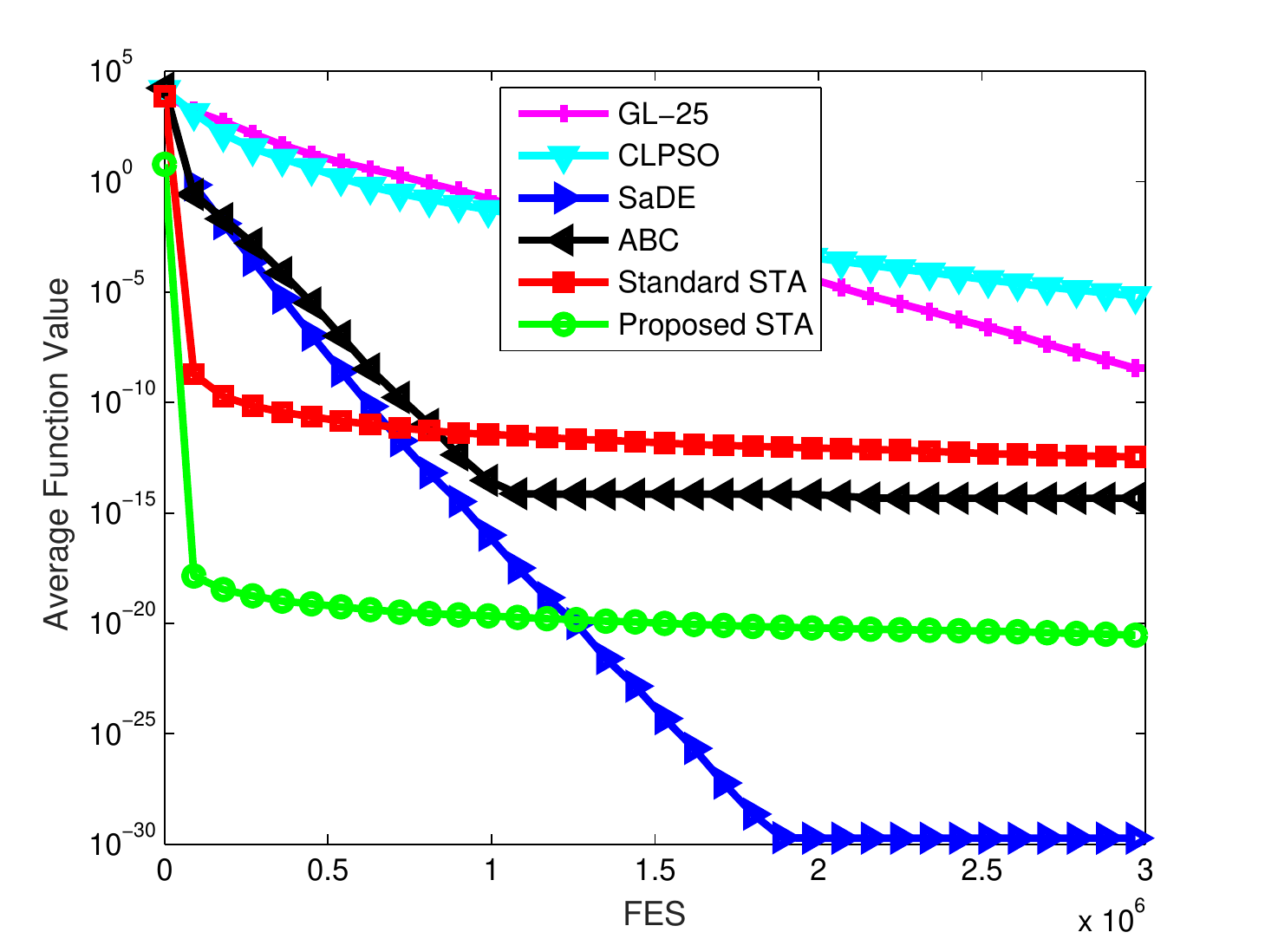}
  \includegraphics[width=6.5cm,height=5cm]{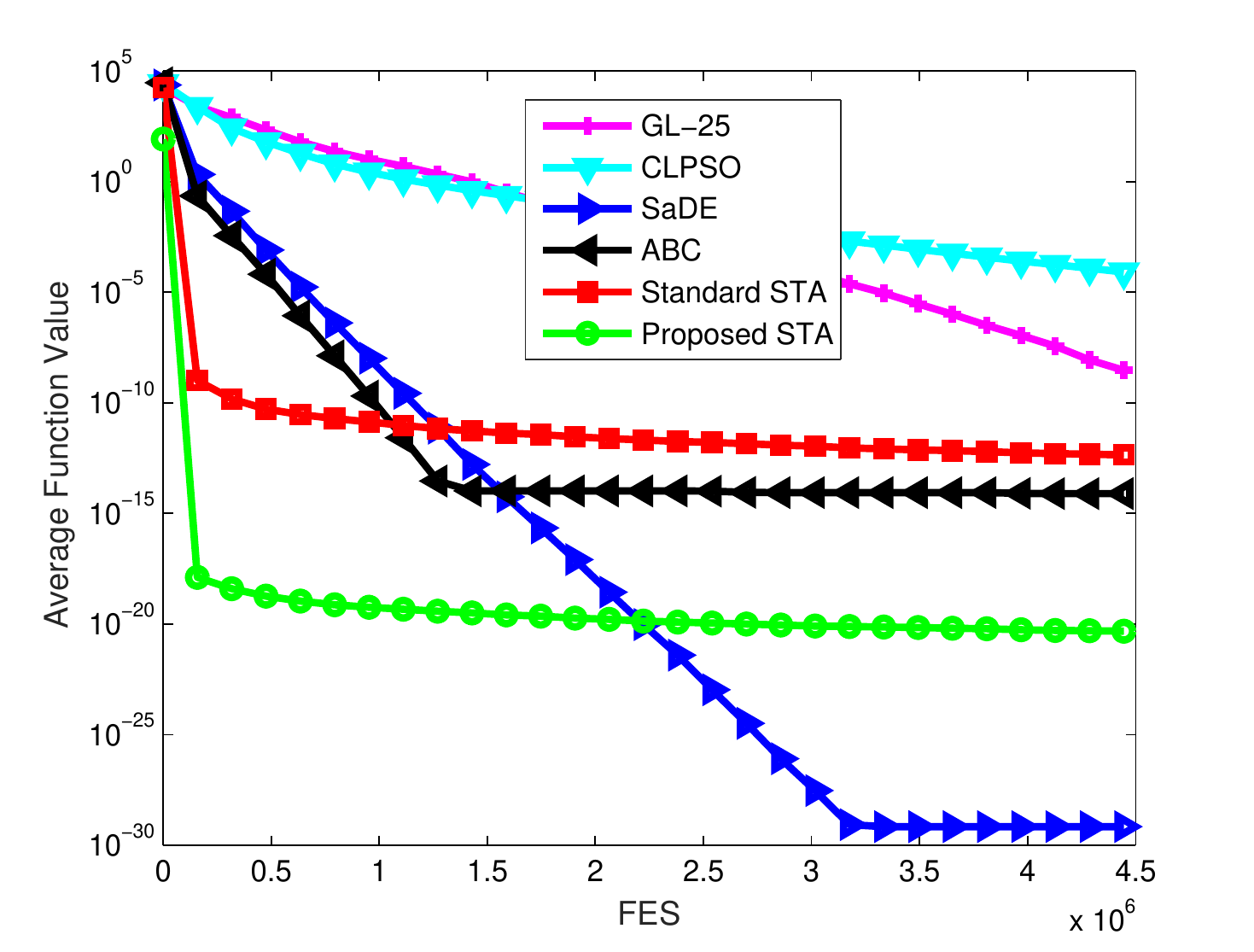}
  \includegraphics[width=6.5cm,height=5cm]{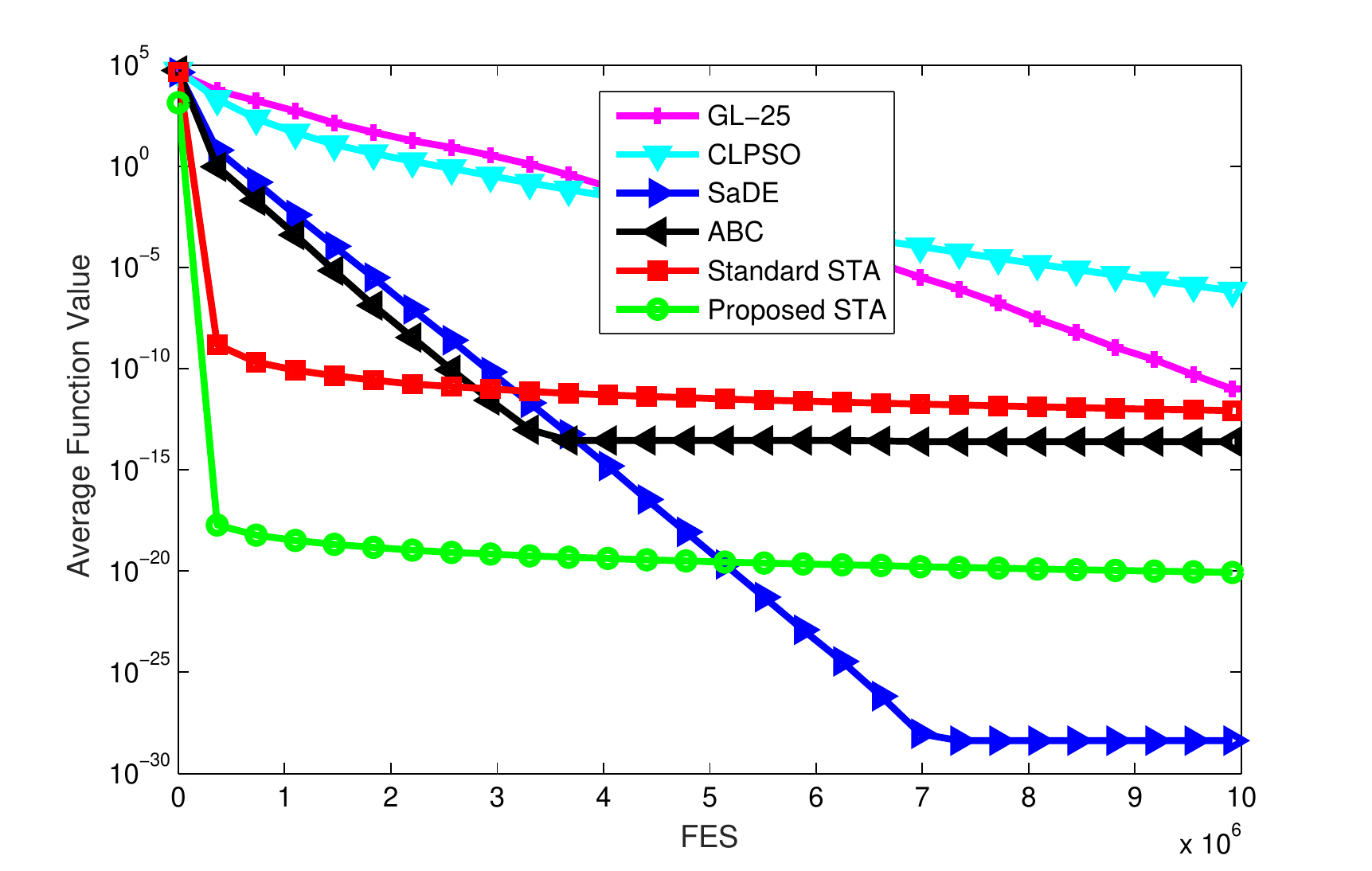}
\caption{The average iterative curves for different metaheuristic methods with respect to the Schwefel 2.4 function}
\label{fig_expansion}       
\end{figure*}

Other metaheuristics are used for comparison, including the GL-25 \cite{garcia2008global},  CLPSO \cite{liang2006comprehensive}, SaDE \cite{qin2009differential}, and ABC \cite{karaboga2007powerful}, with the same parameter settings as in these literatures.
The parameters in the proposed STA are given by experience as follows: \textit{SE} = 30, $T_p = 10$ (additional experiments have testified the validity of these parameter values).
The number of decision variables $n$ of the benchmark functions is set to 20, 30 and 50, and the corresponding
maximum function evaluations is set at 5e4*\textit{n}*log(\textit{n}). A total of 20 independent runs are conducted in the
MATLAB (Version R2010b) software platform on Intel(R) Core(TM) i3-2310M CPU @2.10GHz under Window 7 environment.
The statistic results are given in Table V and some typical instances with respect to elapsed time and iterative curves are illustrated in Figs. 2-4.

From the experimental results, it can be found that the proposed STA is superior to the basic STA among most of these test problems. The global search ability (see the Michalewicz function) and the solution accuracy (see the Schwefel 1.2 and Schwefel 2.4 function) has greatly improved. It can also be comparable to other metaheuristics except for the Michalewicz function and the Schwefel 2.4 function.
However, it should be noted that only mean and standard deviation are given for comparison.
Actually, for the Michalewicz function, the results obtained from the proposed STA hit the known global solution for more than 50\% of the total runs.

\section{Conclusion and future work}
In this study, the optimal parameter selection of operators in continuous STA was considered to improve its search performance. Firstly, a statistical study with four benchmark cases was conducted to investigate how these parameters affect the performance of continuous STA. And several properties are observed from
the statistical study. With the experience gained from the statistical results, then, a new continuous STA with optimal parameters strategy was proposed to accelerate its search process. The proposed STA was successfully applied to other benchmarks. Comparison with other metaheuristics was conducted to demonstrate the effectiveness of the proposed method as well.

It should be noted that the parameter $T_p$ is given by experience that needs further study, and the parameter selection of operators in continuous STA is still a challenging problem, since the proposed optimal parameter selection
strategy can only be considered as a local vision.
From an overall perspective, the parameters set should be taken into consideration as well, and it is not necessarily restricted to one below. Furthermore, it can be found that the STA doesn't work steadily for the Michalewicz function and global search ability should be strengthened further. In our future work, the upper bound of the parameter set will be considered as well, and an adaptive parameter selection
strategy and appropriate utilization of transformation operators can also be alternative choices.

The MATLAB source codes of the standard STA and the proposed STA are available upon request from the corresponding author, or can be downloaded from MATLAB central file exchange, or from X. Zhou's homepage as follows\\ \url{https://www.mathworks.com/matlabcentral/fileexchange/}\\
\url{http://faculty.csu.edu.cn/michael_x_zhou/zh_CN/index.htm}

\section*{Acknowledgment}
This study is supported by the National Natural Science Foundation of
China (Grant No. 61503416, 61533021, 61621062 and 61725306), the Innovation-Driven Plan in Central South University (Grant No. 2018CX12), the 111 Project (Grant No. B17048) and the Hunan Provincial Natural Science
Foundation of China (Grant No. 2018JJ3683).

\ifCLASSOPTIONcaptionsoff
  \newpage
\fi



%
\bibliographystyle{ieeetr}
\bibliography{tzy}

%

\begin{IEEEbiography}[{\includegraphics[width=1in,height=1.25in,clip,keepaspectratio]{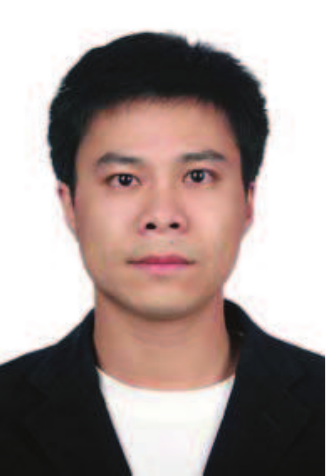}}]{Xiaojun Zhou} received his Bachelor's degree in Automation in 2009 from Central
South University, Changsha, China and received the Ph.D. degree in Applied
Mathematics in 2014 from Federation University Australia.

He is currently an
Associate Professor at Central South University, Changsha, China. His main interests
include modeling, optimization and control of complex industrial process,
optimization theory and algorithms, state transition algorithm, duality theory and their
applications.
\end{IEEEbiography}

\begin{IEEEbiography}[{\includegraphics[width=1in,height=1.25in,clip,keepaspectratio]{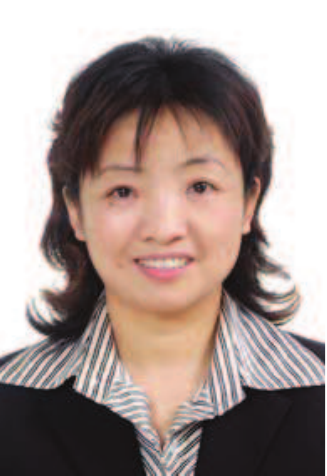}}]{Chunhua Yang} received the M.S.
degree in automatic control engineering and the
Ph.D. degree in control science and engineering
from Central South University, Changsha, China,
in 1988 and 2002, respectively.

From1999 to 2001, she was a Visiting Professor
with the University of Leuven, Leuven, Belgium.
Since 1999, she has been a Full Professor with the
School of Information Science and Engineering,
Central South University. From 2009 to 2010, she
was a Senior Visiting Scholar with the University of Western Ontario, Lon-
don, Canada. Her current research interests include modeling and optimal
control of complex industrial process, fault diagnosis, and intelligent control
system.
\end{IEEEbiography}


\begin{IEEEbiography}[{\includegraphics[width=1in,height=1.25in,clip,keepaspectratio]{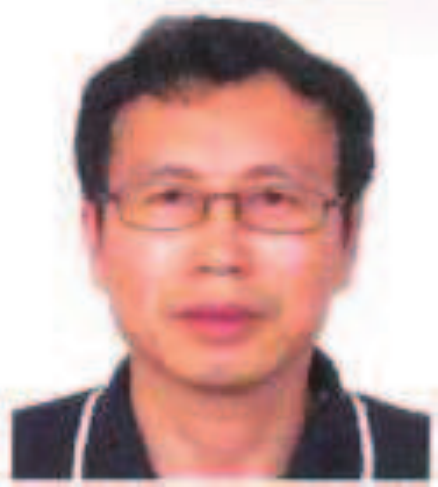}}]{Weihua GUi} received the degree of the B.Eng. in Automatic Control Engineering and
the M.Eng.~ in Control Science and Engineering from Central South University,
Changsha, China, in 1976 and 1981, respectively.

From 1986 to 1988, he was a
visiting scholar at Universitat-GH-Duisburg, Germany. He is a member of the
Chinese Academy of Engineering and has been a full professor in the School of
Information Science and Engineering, Central South University, Changsha, China,
since 1991. His main research interests are in modeling and optimal control of
complex industrial process, distributed robust control, and fault diagnoses.
\end{IEEEbiography}




\end{document}